\documentclass[12pt, leqno]{amsart}

\usepackage{graphicx}
\usepackage{amssymb,amsmath,amsthm,amscd}
\usepackage{mathrsfs}
\usepackage{enumerate}
\usepackage{enumitem}
\usepackage{verbatim}
\usepackage[usenames,dvipsnames]{color}
\usepackage[colorlinks=true, pdfstartview=FitV,
linkcolor=blue,citecolor=blue,urlcolor=blue]{hyperref}
\usepackage[all]{xy}
\usepackage{tikz}
\usetikzlibrary{shapes, arrows}
\usepackage{tikz-cd}
\usetikzlibrary{matrix, arrows}
\usepackage{stmaryrd} 
\usepackage{dsfont}
\usepackage{bm} 
\usepackage{scalerel}[2016/12/29] 
\usepackage{amsthm}
\usepackage{adjustbox}
\usepackage{float}
\usepackage{array}
\usepackage{multirow} 
\usepackage{caption} 
\allowdisplaybreaks[4]

\DeclareMathSymbol{\shortminus}{\mathbin}{AMSa}{"39}

\newcommand{\nc}{\newcommand}

\numberwithin{equation}{section}

\usepackage[colorinlistoftodos]{todonotes}
\usepackage{mathtools}
\mathtoolsset{showonlyrefs,showmanualtags}
\usepackage{empheq}

\theoremstyle{plain}
\newtheorem{lem}{Lemma}[section]
\newtheorem{pro}[lem]{Proposition}
\newtheorem{thm}[lem]{Theorem}
\newtheorem{cor}[lem]{Corollary}
\newtheorem{defi}[lem]{Definition}

\newcommand{\Pro}{\begin{pro}}
	\newcommand{\enpro}{\end{pro}}
\newcommand{\Lem}{\begin{lem}}
	\newcommand{\enlem}{\end{lem}}
\newcommand{\Thm}{\begin{thm}}
	\newcommand{\enthm}{\end{thm}}
\newcommand{\Cor}{\begin{cor}}
	\newcommand{\encor}{\end{cor}}
\newcommand{\Defi}{\begin{defi}}
	\newcommand{\enDefi}{\end{defi}}
\newcommand{\Proof}{\begin{proof}}
	\newcommand{\enproof}{\end{proof}}

\theoremstyle{definition} 
\newtheorem{rem}[lem]{Remark}
\newtheorem{exam}[lem]{Example}
\newtheorem{Convention}[lem]{Convention}

\newcommand{\Conv}{\begin{Convention}}
	\newcommand{\enconv}{\end{Convention}}
\nc{\Rem}{\begin{rem}}
	\nc{\enrem}{\end{rem}}

\newcommand{\arxiv}[1]{\href{http://arxiv.org/abs/#1}{\tt arXiv:\nolinkurl{#1}}}

\newcommand{\monoto}{\hookrightarrow}
\nc{\epito}{\twoheadrightarrow}
\newcommand{\isoto}[1][]{\mathop{\xrightarrow%
		[{\raisebox{.3ex}[0ex][.3ex]{$\scriptstyle{#1}$}}]%
		{{\raisebox{-.6ex}[0ex][-.6ex]{$\mspace{2mu}\sim\mspace{2mu}$}}}}}

\nc{\rmkend}{\hfill$\triangledown$}
\nc{\defend}{\hfill$\triangle$}

\nc{\ccc}{\mathfrak{c}}
\nc{\CCC}{\mathfrak{C}}
\nc{\Ck}{\mathfrak{C}}

\nc{\kor}{\mathbb{C}}
\nc{\indx}{\mathbb{I}}

\nc{\CC}{C}
\nc{\cc}{c}
\nc{\sss}{s}
\nc{\ck}{\mathfrak{c}}

\nc{\Bg}{B}
\nc{\Ag}{A}
\nc{\As}{\pmb{\Ag}(z)}
\nc{\Asmone}{\pmb{\Ag}^{(-1)}(z)}
\nc{\Aps}{\pmb{\Ag}_+(z)}
\nc{\Apsone}{\pmb{\Ag}^{(1)}_+(z)}
\nc{\Apsmone}{\pmb{\Ag}^{(-1)}_+(z)}
\nc{\Ams}{\pmb{\Ag}_-(z)}
\nc{\Hg}{H}
\nc{\Thg}{\Theta}
\nc{\Thgs}{\pmb{\Theta}}
\nc{\Thgsr}{\pmb{\grave{\Theta}}}
\nc{\bTh}{\pmb{\grave{\Theta}}(z)}
\nc{\Apsvar}[1]{\pmb{\Ag}^{(#1)}_+(z)} 

\nc{\smin}{{\shortminus}}

\nc{\Ui}{\widetilde{\mathbf{U}}^\imath}
\nc{\DrUi}{{}^{\mathrm{Dr}}\Ui}
\nc{\fext}[2]{{#1}[\negthinspace[#2]\negthinspace]}

\nc{\gTh}{\grave{\Theta}}

\DeclareRobustCommand{\qbinom}{\genfrac{[}{]}{0pt}{}}

\nc{\DD}{\pmb{D}}
\nc{\KK}{\mathbb{K}}

\nc{\xgp}{x^+}
\nc{\Xgp}{\pmb{x}^+}
\nc{\xgm}{x^-}
\nc{\Xgm}{\pmb{x}^-}
\nc{\xgpm}{x^{\pm}}
\nc{\psig}{\phi^+}
\nc{\phig}{\phi^-}
\nc{\phipm}{\phi^\pm}
\nc{\Psig}{\pmb{\phi^+}}
\nc{\Phig}{\pmb{\phi^-}}
\nc{\Phipm}{\pmb{\phi^\pm}}
\nc{\hg}{h}

\nc{\bY}{\mathbf{Y}}
\nc{\bA}{\mathbf{A}}
\nc{\bT}{\mathbf{T}}

\nc{\Eg}{e}
\nc{\Kg}{K}

\nc{\Rep}{\on{Rep}}

\nc{\qq}{(q-q^{-1})^{-1}}
\nc{\factor}{\Omega}

\nc{\chring}{\Z[Y_a^{\pm 1}]_{a \in \C^\times}}
\nc{\chmap}{\chi_q}
\nc{\chmod}{\Z[\mathbf{Y}_a^{\pm 1}]_{a \in \C^\times}}
\nc{\ourchmap}{\boldsymbol{\chi}_q}
\nc{\Yring}{\mathcal{Y}} 
\nc{\Ymod}{\pmb{\mathcal{Y}}}
\nc{\Yg}{\mathbf{Y}_{i,a}}
\nc{\Ys}{\mathbf{Y}_{i,\mathbf{s}}}
\nc{\Kmat}{\mathcal{K}^0} 

\nc{\Tr}{\on{Tr}}
\nc{\id}{\on{id}}

\nc{\ourR}{\mathbf{R}}
\nc{\ourQ}{\mathbf{Q}}
\nc{\ourY}{\mathbf{Y}}

\nc{\Oq}{\mathcal{O}_q(\widehat{\g})}
\nc{\Oqi}{\mathcal{O}_q^{[i]}}
\nc{\Oqc}{\mathcal{O}_q^{\mathbf{c}}(\widehat{\g})}
\nc{\car}{\mathcal{H}}
\nc{\qaa}{U_q(\widehat{\mathfrak{g}})}
\nc{\qla}{U_q(L\mathfrak{g})}
\nc{\uqla}{\widetilde{U}_q(L\mathfrak{g})}
\nc{\drqaa}{{}^{\mathrm{Dr}}\qaa}
\nc{\uqsl}{U_q(\widehat{\mathfrak{sl}}_2)}
\nc{\uLsl}{U_qL\mathfrak{sl}_2}
\nc{\Serre}{\mathsf{Serre}}
\nc{\Sym}{\on{Sym}}
\nc{\UXp}{UX_+} 

\nc{\Tbr}{\mathbf{T}}
\nc{\degdr}{\on{deg}^{\mathrm{Dr}}}
\nc{\Uq}{\mathbf{U}}
\nc{\Uu}{\widetilde{\mathbf{U}}}
\nc{\gaf}{\widehat{\mathfrak{g}}}

\nc{\ada}{\ad_{\bhg_{-1}}}
\nc{\adb}{\ad_{\bhg_{1}}}
\nc{\adc}{\ad_{\bHg_{1}^{(2)}}}

\nc{\tAg}{\widetilde{A}}
\nc{\hAg}{\widehat{A}}

\nc{\bHg}{\overline{\Hg}}
\nc{\bhg}{\overline{\hg}}
\nc{\tHg}{\widetilde{H}}

\nc {\DrOq}{{}^{\mathrm{Dr}}\Oq}
\nc{\gup}[1]{^{(#1)}}
\nc{\gupp}{^{(1),+}}
\nc{\mi}{^{-1}}
\nc{\adh}{\operatorname{ad}_{\bar{\hg}_{-1}}}
\nc{\adhp}{\operatorname{ad}_{\bar{\hg}_{1}}}
\nc{\Omg}{\Omega^{-1}}
\nc{\Htwo}{\overline{H}_1^{(2)}}
\nc{\ad}{\operatorname{ad}}

\nc{\sspan}{\on{span}}

\newcommand{\commentout}[1]{}

\newcommand{\on}{\operatorname}
\nc{\be}{\begin{enumerate}}
	\nc{\ee}{\end{enumerate}}
\newcommand{\eq}{\begin{equation}}
	\newcommand{\eneq}{\end{equation}}
\nc{\bc}{\begin{cases}}
	\nc{\ec}{\end{cases}}
\newcommand{\eqn}{\begin{eqnarray*}}
	\newcommand{\eneqn}{\end{eqnarray*}}
\newcommand{\ba}{\begin{array}}
	\newcommand{\ea}{\end{array}}

\newcommand{\C}{{\mathbb C}}

\newcommand{\Z}{{\mathbb Z}}
 
\newcommand{\R}{{\mathbb R}} 

\newcommand{\g}{{\mathfrak{g}}}
\nc{\Ad}{\operatorname{Ad}}

\nc{\gr}{\on{gr}}

\newcommand{\End}{\operatorname{End}}
\nc{\Aut}{\operatorname{Aut}}

\nc{\coker}{\operatorname{coker}}

\nc{\Img}{\on{Im}}
\nc{\res}{\on{res}}

\nc{\modv}[1]{{#1}\operatorname{-mod}}

\nc{\bl}{\bigl(}
\nc{\br}{\bigr)}

\newlength{\mylength}
\DeclareRobustCommand{\SkipTocEntry}[5]{}

\AtBeginDocument{%
   \def\MR#1{}
}

\title[Compatibility of Drinfeld presentations for QSP]
{Compatibility of Drinfeld presentations for split affine Kac--Moody quantum symmetric pairs}

\author[J.-R. Li]{Jian-Rong Li}
\address{Faculty of Mathematics, University of Vienna, Oskar Morgenstern Platz 1, 1090 Vienna, Austria}
\email{\href{mailto:lijr07@gmail.com}{lijr07@gmail.com}}

\author[T. Prze\'{z}dziecki]{Tomasz Prze\'{z}dziecki}
\address{School of Mathematics, University of Edinburgh, Peter Guthrie Tait Rd, Edinburgh, EH9 3FD, United Kingdom, OrciD: 0000-0001-9700-1007}
\email{\href{mailto:tprzezdz@exseed.ed.ac.uk}{tprzezdz@exseed.ed.ac.uk}}



\thanks{The first author was supported by the Austrian Science Fund (FWF): P-34602, Grant DOI: 10.55776/P34602, and PAT 9039323, Grant-DOI 10.55776/PAT9039323. The second author was supported by the EPSRC grant No.\ EP/W022834/1 \emph{Kac--Moody quantum symmetric pairs, KLR algebras and generalized Schur--Weyl duality}.}

\begin{document}

\begin{abstract}
Let $(\Uq, \Uq^\imath)$ be a split affine quantum symmetric pair of type $\mathsf{B}_n^{(1)}, \mathsf{C}_n^{(1)}$ or~$\mathsf{D}_n^{(1)}$. We prove factorization and coproduct formulae for the Drinfeld--Cartan operators $\Theta_i(z)$ in the Lu--Wang Drinfeld-type presentation, generalizing the type $\mathsf{A}_n^{(1)}$ result from \cite{Przez-23}. As an application, we show that a boundary analogue of the $q$-character map, defined via the spectra of these operators, is compatible with the usual $q$-character map. As an auxiliary result, we also produce explicit reduced expressions for the fundamental weights in the extended affine Weyl groups of classical types. 
\end{abstract}

\maketitle

\setcounter{tocdepth}{1}
\tableofcontents

\section{Introduction}
\nc{\Uh}{U_q(\widetilde{\mathfrak{h}})}
\nc{\Uhz}{\fext{U_q(\widetilde{\mathfrak{h}})}{z}}
\nc{\ichmap}{\chi_q^\imath}

Quantum affine algebras admit three distinct presentations: the `Drinfeld--Jimbo', `new Drinfeld' and `RTT' realizations. The first can be seen as a quantization of the usual Serre-type presentation of a Kac-Moody Lie algebra $\widehat{\g}$, while the second is a quantization of the central extension presentation of $\widehat{\g}$. The interplay between these different realizations was studied and described precisely in \cite{Ding-Fr, FrMukHopf}. One of the most interesting features of the new Drinfeld presentation is that it exhibits a large, infinitely generated, commutative subalgebra of $U_q(L\g)$. The spectra of the generators of this subalgebra, often called Drinfeld--Cartan operators $\phi^\pm_{i,m}$, play a key role in the classification of finite dimensional representation of $U_q(L\g)$ via Drinfeld polynomials \cite{chari-pressley-qaa, chari-pressley-94}, and $q$-character theory \cite{FrenRes, FrenMuk-comb}. 

Recently, Lu and Wang \cite{lu-wang-21}, building on the work of Baseilhac and Kolb \cite{bas-kol-20} in rank one (see also \cite{ZhangDr, LWZ-quasi}), have constructed a Drinfeld-type presentation for split and certain quasi-split affine quantum symmetric pair coideal subalgebras. The Lu--Wang presentation also exhibits a large commutative subalgebra, generated by $\Theta_{i,m}$. 
This very significant development raises many questions. Here we focus on the following three. \emph{Firstly}, what is the relationship between the generators $\Theta_{i,m}$ and the usual Drinfeld--Cartan operators $\phi^\pm_{i,m}$? \emph{Secondly}, how to describe the coproduct of $\Theta_{i,m}$? \emph{Thirdly}, can the Lu--Wang presentation be used to define a `boundary' analogue of $q$-characters? 

The first two questions were answered by the second author \cite{Przez-23} for split affine quantum symmetric pairs of type $\mathsf{A}_n^{(1)}$. In particular, it was shown that the generating series $\Thgsr_i(z)$ have the following \emph{factorization property}: 
\eq \label{intro 1}
\Thgsr_i(z) \equiv \pmb{\phi}_i^-(z\mi)\pmb{\phi}_i^+(C z) \quad \mod \fext{U_q(L\g)_+}{z}, 
\eneq
where $U_q(L\g)_+$ is the Drinfeld positive half of $U_q(L\g)$. Moreover, the series $\Thgsr_i(z)$ are \emph{approximately group-like}, in the sense that 
\eq \label{intro 2}
\Delta(\Thgsr_i(z)) \equiv \Thgsr_i(z) \otimes \Thgsr_i(z) \quad  \mod \fext{U_q(L\g)  \otimes U_q(L\g)_{+}}{z}. 
\eneq
The main result of the present paper is a generalization of the two results above to split affine quantum symmetric pairs of types $\mathsf{B}_n^{(1)}, \mathsf{C}_n^{(1)}$ and $\mathsf{D}_n^{(1)}$. As an application, we also propose a definition of boundary $q$-characters and prove that they are, in a suitable sense, compatible with the usual $q$-characters.

\subsection{Proof strategy} 

Our proof is based on an algebraic and computational approach, using the Drinfeld--Jimbo and new Drinfeld presentations only. In particular, we do not use any information available through the RTT presentation. Since the overall proof is quite technical, let us summarize the main steps below. 

\emph{Step 1: compute explicit reduced expressions for the fundamental weights $\omega_i$.} This is carried out in \S \ref{sec-comb}. Whenever $\alpha_i$ has multiplicity one in the highest root, the expressions can be derived from a formula in \cite{bourbaki}. In the remaining cases, we consider the extended affine Weyl group as a group of affine transformations of $\R^n$ and compute that our expressions act as the correct affine shifts. To show that they are reduced, we calculate their lengths and compare them with a formula from \cite{IM65}. 

\emph{Step 2: express $\Tbr_{\omega'_i}(B_i)$ as an explicit polynomial in the variables $B_j$.} As in the usual Drinfeld presentation, the generators in the Lu--Wang presentation are constructed using a braid group action (which differs from the ordinary Lusztig action). More precisely, they are constructed by repeatedly applying the braid group operators $\Tbr_{\omega_i}$ to the usual Kolb--Letzter generators of the coideal subalgebra. Using the explicit formulae for $\omega_i$ from Step $1$, we are able to express $\Tbr_{\omega'_i}(B_i)$ by means of certain recursively defined polynomials in the variables $B_j$. This is carried out in \S \ref{sec: weak com}, based on the auxiliary calculations from \S \ref{sec: root combin}. 

\emph{Step 3: show that $\Tbr_{\omega'_i}(B_i) \equiv T_{\omega'_i}(B_i)$.} More precisely, we show that the Lu--Wang braid group operator $\Tbr_{\omega'_i}$ coincides with the usual Lusztig braid group operator $T_{\omega'_i}$ on $B_i$, modulo a certain subalgebra of $U_q(L\g)_+$. This is achieved in \S \ref{sec: weak com} (Theorem \ref{thm: weak comp}) in the following way. Recall that $B_j$ is a linear combination of $F_j$ and $E_jK_j\mi$. Writing $\Tbr_{\omega'_i}(B_i)$ as a polynomial as in Step $2$, and substituting for each occurrence of $B_j$ either $F_j$ or $E_jK_j\mi$, we can express $\Tbr_{\omega'_i}(B_i)$ as a sum of homogeneous terms in the Drinfeld gradation. The main problem is to show that, except for $T_{\omega'_i}(E_iK_i\mi)$, the terms which are not in $U_q(L\g)_+$ vanish. This is solved using the theory of good polynomials developed in \S \ref{sec: good polynomials}, with the aid of auxiliary calculations from \S \ref{sec: aux cal brgr}.  

\emph{Step 4: deduce the factorization and coproduct formulae.} The main idea is that, modulo $U_q(L\g)_+$, we can reduce the problem to rank one subalgebras, and apply the results from \cite{Przez-23}. We show that, to perform this reduction, we only need to verify a single commutation condition. This criterion is then checked explicitly in type $\mathsf{D}_n^{(1)}$. The same methods carry over to types $\mathsf{B}_n^{(1)}$ and $\mathsf{C}_n^{(1)}$. This is achieved in \S \ref{sec: strong compat}. 

We remark that the overall proof strategy is similar to that in \cite[\S 9]{Przez-23}. However, each step of the proof is significantly harder than in type $\mathsf{A}$ case, mainly due to the higher complexity of reduced expressions for fundamental weights in other classical types. Another innovation is the concept of `good polynomials', which allows us to systematically handle different types, and is also likely to find application in, e.g., the setting of quasi-split affine quantum symmetric pairs. 

\subsection{Reduced expressions for fundamental weights} 

In Step $1$ above, we compute explicit reduced expressions for the fundamental weights in the extended affine Weyl groups of classical types. Somewhat surprisingly, to our knowledge, such explicit expressions have not appeared in the literature before\footnote{A recursive formula can be found in \cite{Lus83}.}. As this result is interesting in its own right, we collect the expressions in a list below for quick reference, see Table \ref{tab:expresssion of fundamental weights of the extended affine Weyl group}.

\begin{table}[h]
\centering
\begin{tabular}{|c|c|c|}
\hline
Type & Fundamental weight $\omega_i$ & Range of $i$ \\
\hline
$\mathsf{A}_n^{(1)}$ & $\pi^i [n {\shortminus} i{+}1 ,n] \cdots[2,i{+}1][1,i]$ & $1 \leq i \leq n$ \\
\hline
\multirow{2}{*}{$\mathsf{B}_n^{(1)}$} & $\pi_1 [1,n] \big(s_0[2,n][1,n]\big)^{\frac{i {\shortminus} 1}{2}} [{n{\shortminus}i},{n{\shortminus}1}] \cdots [1,i]$ & $1 \leq i \leq n$, $i$ is odd \\
& $\big(s_0[2,n][1,n]\big)^{\frac{i}{2}}  [{n{\shortminus}i},{n{\shortminus}1}] \cdots [1,i]$ & $2 \leq i \leq n$, $i$ is even \\
\hline
\multirow{2}{*}{$\mathsf{C}_n^{(1)}$} & $\pi_n s_n[n {\shortminus} 1,n] \cdots [1,n]$ & $i=n$ \\
& $(s_0 [1,n] )^i [n {\shortminus} i,n {\shortminus} 1] \cdots [2,i{+}1][1,i]$ & $i<n$ \\
\hline
\multirow{6}{*}{$\mathsf{D}_n^{(1)}$} & $\big(s_0 [2,n{\shortminus}1][1,n{\shortminus}2]s_n\big)^{\frac{i}{2}} [{n{\shortminus}i},{n{\shortminus}1}] \cdots [1,i] $ & $2 \leq i \leq n {\shortminus} 2$, $i$ is even \\
& $\pi_{n{\shortminus}1} r_{n{\shortminus}2}r_{n{\shortminus}4} \cdots r_3 s_n [1,{n{\shortminus}1}]$ & $i=n {\shortminus} 1$, $i$ is even \\
& $\pi_{n} r_{n{\shortminus}2} r_{n{\shortminus}4} \cdots r_2 s_n$ & $i=n$, $i$ is even \\
& $\pi_1[1,n{\shortminus}2]s_n\big(s_0[2,n{\shortminus}1][1,n{\shortminus}2]s_n\big)^{\frac{i{\shortminus}1}{2}}[{n{\shortminus}i},{n{\shortminus}1}] \cdots [1,i]$ & $1 \leq i \leq n {\shortminus} 2$, $i$ is odd \\
& $\pi_{n{\shortminus}1} s_{n{\shortminus}1}s_{n{\shortminus}2} r_{n{\shortminus}3} r_{n{\shortminus}5} \cdots r_3 s_n [1,{n{\shortminus}1}]$ & $i=n{\shortminus}1$, $i$ is odd \\
& $\pi_{n} s_{n{\shortminus}1}s_{n{\shortminus}2}r_{n{\shortminus}3} r_{n{\shortminus}5} \cdots r_2 s_n$ & $i=n$, $i$ is odd \\
\hline
\end{tabular}

\caption{Reduced expressions for fundamental weights $\omega_i$ in extended affine Weyl groups of classical type. Here $[k,l] = s_k s_{k+1} \cdots s_{l}$ and $r_m = s_n [m,n{\shortminus}1][m{\shortminus}1,n{\shortminus}2]$.}
\label{tab:expresssion of fundamental weights of the extended affine Weyl group}
\end{table}

\subsection{Application to $q$-characters}

The notion of $q$-characters can be defined in at least three equivalent ways: via the universal $R$-matrix, via the spectrum of Drinfeld--Cartan operators, or via Nakajima's quiver varieties. The quantum symmetric pair analogues of these approaches are objects of intense study, see, e.g., \cite{AppelVlaar2, AppelVlaar, AppelVlaar3} for the latest developments on the universal $K$-matrix, and \cite{YiqiangLi} for a quiver variety approach to symmetric pairs. Nevertheless, it appears that the geometric and integrable systems methods are not yet advanced enough to yield a satisfactory theory of $q$-characters for quantum symmetric pairs. The main obstacle in the latter case is the unavailability of a Khoroshkin-Tolstoy-Levendorsky-Soibelman-Stukopkin-Damiani-type factorization of the (affine)  universal $K$-matrix. 

In light of the aforementioned difficulties, we propose to define boundary $q$-characters directly via the Lu--Wang presentation instead. More precisely, we consider the generalized eigenspace decomposition of a finite dimensional representation with respect to the action of the operators $\Theta_{i,m}$, and let boundary $q$-characters encode the multiplicities of such eigenspaces. This is equivalent to taking the trace of a certain operator in the completion $\fext{\Uq^\imath \otimes \Uh}{z}$, yielding a map $\ichmap \colon \Rep \Uq^\imath \to \Uhz$, where 
$U_q(\widetilde{\mathfrak{h}}) = \langle k_i, h_{i,r} \mid i \in \indx_0, r \leq 0 \rangle$ is a `half' of the Drinfeld--Cartan subalgebra of $U_q(L\g)$. We apply \eqref{intro 1}--\eqref{intro 2} to show that $\ichmap$ is compatible with the usual $q$-character map (Corollary \ref{cor: comm diagram qchar actions}). More precisely, we show that the diagram 
\[
\begin{tikzcd}[ row sep = 0.2cm]
\Rep U_q(L\g) \arrow[r, "\chmap"] & \Z[Y_{i,a}^{\pm 1}]_{i \in \indx_0, a \in \C^{\times}}  \\
 \curvearrowright  & \curvearrowright  \\
\Rep \Uq^\imath \arrow[r, "\ichmap"] & \Uhz 
\end{tikzcd}
\]
commutes if $\Uhz$ is endowed with an appropriate `twisting' action of $\Z[Y_{i,a}^{\pm 1}]_{i \in \indx_0, a \in \C^{\times}}$. 
In particular, this result yields an easy way to compute the boundary $q$-characters of restriction representations. 
For first results in the case of more general representations, we refer the reader to \cite{LP25}.

\addtocontents{toc}{\SkipTocEntry}
	
\section*{Acknowledgements} 
We would like to thank Andrea Appel, Pascal Baseilhac and Bart Vlaar for insightful discussions and comments. 

\section{Preliminaries}

We work over the field of complex numbers and assume that $q \in \C^\times$ is not a root of unity throughout. 

\subsection{Quantum affine algebras} 

Let $\indx_0 = \{ 1,\cdots, n \}$ and $\indx = \indx_0 \cup \{0\}$. Let $\g$ be a simple Lie algebra with Cartan matrix $(a_{ij})_{i,j \in \indx_0}$, and $\widehat{\g}$ the corresponding untwisted affine Lie algebra with affine Cartan matrix $(a_{ij})_{i,j \in \indx}$. Let $d_i$ be relatively prime positive integers such that $(d_ia_{ji})_{i,j\in\indx_0}$ is a symmetric matrix, and set $q_i = q^{d_i}$. We use Bourbaki's conventions \cite{bourbaki} for the Cartan matrix in non-simply laced types. i.e., 
\begin{align*}
\left(\begin{array}{rrrrrrr}
2 & -1 & \hdots  & 0 & 0 \\
-1 & 2 & \hdots   & 0 & 0 \\
\vdots & \vdots &  \ddots  & \vdots & \vdots \\
0 & 0 & \hdots   & 2 & -2 \\
0 & 0 & \hdots   & -1 & 2
\end{array}\right) \quad , \quad 
\left(\begin{array}{rrrrrrr}
2 & -1 & \hdots  & 0 & 0 \\
-1 & 2 & \hdots   & 0 & 0 \\
\vdots & \vdots &  \ddots  & \vdots & \vdots \\
0 & 0 & \hdots   & 2 & -1 \\
0 & 0 & \hdots   & -2 & 2
\end{array}\right) 
\end{align*}
are the Cartan matrices of types  $\mathsf{B}_n$ and  $\mathsf{C}_n$, respectively. 

We use standard notations and conventions regarding root systems, weight lattices, Weyl groups, etc., as in, e.g., \cite[\S 3.1]{lu-wang-21}. In particular, we 
let $\alpha_i$ ($i \in \indx$) denote the simple roots of $\widehat{\g}$; let $\theta$ denote the highest root of $\g$, and $\delta$ the basic imaginary root; let $P$ and $Q$ denote the weight and root  lattices of $\g$, respectively; and let $\omega_i \in P$ ($i \in \indx_0$) be the fundamental weights of~$\g$. 
Note that, according to Bourbaki's conventions, $s_{\alpha_i}(\alpha_j) = \alpha_j - a_{ji} \alpha_i$. 

The \emph{quantum affine algebra} $\qaa$ is the algebra with generators $\Eg_i^\pm, \Kg_i^{\pm 1}$ $(i \in \indx)$ and relations: 
\begin{align}
\Kg_i\Kg_i^{-1} &= \Kg_i^{-1}\Kg_i = 1, \\ 
\Kg_i\Kg_j &= \Kg_j \Kg_i, \\
\Kg_i\Eg_j^{\pm} &= q_i^{\pm a_{ji}} \Eg_j^\pm\Kg_i,    \label{eq: Ke eK rel} \\
[\Eg_i^{+}, \Eg_j^{-}] &= \delta_{ij} \frac{\Kg_i - \Kg_i^{-1}}{q_i - q_i^{-1}},   \label{eq: Ei Ej rel} \\
\Serre_{ij}(\Eg_i^{\pm}, \Eg_j^\pm) &= 0 \qquad (i \neq j),
\end{align}
where
\[
\Serre_{ij}(x,y) = \sum_{r=0}^{1-a_{ji}} (-1)^r \qbinom{1-a_{ji}}{r}_{q_i} x^{1-a_{ji}-r}yx^r. 
\] 
We will also sometimes abbreviate
\[
E_i = e_i^+, \qquad F_i = e_i^-.
\]
We use the standard notation for divided powers, i.e., $(e^{\pm}_i)^{(r)} = (e^{\pm}_i)^{r}/[r]_{q_i}!$. 
Given $\mu = \sum_{i \in \indx} c_i \alpha_i \in \Z \indx = \bigoplus_{i \in \indx} \Z \alpha_i$, set 
\eq \label{eq: Kmu defi}
K_\mu = \textstyle \prod_{i \in \indx} K_i^{c_i}, \qquad K_{\delta} = K_0 K_{\theta}. 
\eneq

The algebra $\qaa$ is a Hopf algebra, with the coproduct
\eq \label{eq: coprod on U}
\Delta(\Eg_i^+) = \Eg_i^+ \otimes 1 + \Kg_i \otimes \Eg_i,  \quad \Delta(\Eg_i^-) = \Eg_i^- \otimes \Kg_i^{-1} + 1 \otimes \Eg_i^-, \quad \Delta(\Kg_i^{\pm 1}) = \Kg_i^{\pm 1} \otimes \Kg_i^{\pm 1}. 
\eneq
The counit is given by 
\[
\varepsilon(\Eg_i^{\pm}) = 0, \quad \varepsilon(\Kg_i^{\pm 1}) = 1. 
\]

The \emph{quantum loop algebra} $\qla$ is the quotient of $\qaa$ by the ideal generated by the central element $\Kg_\delta - 1$. 
Let us recall the ``new" Drinfeld presentation of the quantum loop algebra $\qla$. By \cite{drinfeld-dp, beck-94}, $\qla$ is isomorphic to the algebra generated by $\xgpm_{i,k}, \hg_{i,l}, \Kg^{\pm 1}_{i}$, where $k \in \Z$, $l \in \Z - \{0\}$ and $i \in \indx_0$, subject to the following relations:
\begin{align}
\Kg_i\Kg_i^{-1} =& \ \Kg_i^{-1}\Kg_i = 1, \\
\Kg_i \Kg_j =& \ \Kg_j \Kg_i, \\ 
[\hg_{i,k}, \hg_{j,l}] =& \ 0, \\
\Kg_i\hg_{j,k} =& \ \hg_{j,k}\Kg_i, \\ 
\Kg_i \xgpm_{j,k} =& \ q_i^{\pm a_{ji}} \xgpm_{j,k}\Kg_i, \label{eq: Kx rel} \\
[\hg_{i,k},\xgpm_{j,l}] =& \pm \textstyle \frac{[k\cdot a_{ji}]_{q_i}}{k}\xgpm_{j,k+l} \label{eq: hx rel},\\
\xgpm_{i,k+1}\xgpm_{j,l} - q_i^{\pm a_{ji}} \xgpm_{j,l}\xgpm_{i,k+1} =& \ q_i^{\pm a_{ji}}\xgpm_{i,k}\xgpm_{j,l+1} - \xgpm_{j,l+1}\xgpm_{i,k}, \label{eq: sl2 4x rel} \\
[\xgp_{i,k}, \xgm_{j,l}] =& \ \delta_{ij} \textstyle\frac{1}{q_i-q_i^{-1}}(\psig_{i, k+l} - \phig_{i, k+l}), \label{eq: g x+-}
\end{align}
\vspace{-0.25cm}
\[
\Sym_{k_1,\cdots,k_r} \sum_{s=0}^r (-1)^s \qbinom{r}{s}_{q_i} \xgpm_{i,k_1} \cdots \xgpm_{i,k_s}\xgpm_{j,l}\xgpm_{i,k_{s+1}} \cdots \xgpm_{i,k_r} = 0  \quad \mbox{for } r=1-a_{ji} \quad (i \neq j), 
\]
where $\Sym_{k_1,\cdots,k_r} $ denotes symmetrization with respect to the indices $k_1, \cdots, k_r$ and 
\[
\pmb{\phi}^{\pm}_i(z) = \sum_{k=0}^\infty \phipm_{i,\pm k} z^{\pm k} = \Kg_i^{\pm 1} \exp\left( \pm (q-q^{-1}) \sum_{k=1}^\infty \hg_{i,\pm k} z^{\pm k} \right). 
\]

Let 
\[
\uqla = \qla \otimes \C[\KK_i^{\pm1} \mid i \in \indx].  
\]
It can be regarded as a version of the Drinfeld double of $U_q(\widehat{\mathfrak{n}}_+)$. 
Given $\mu \in \Z \indx$, we define $\KK_\mu$ in the same way as in \eqref{eq: Kmu defi}. 
We consider $\uqla$ as a $Q$-graded algebra with 
\[
\degdr \xgpm_{i,k} = \pm \alpha_i, \quad \degdr h_{i,k} = \degdr \KK_i = 0 \qquad (i \in \indx_0). 
\] 
In the following, we will use the abbreviations
\[
\Uq = \qla, \qquad \Uu = \uqla. 
\] 
Let $\Uu_+$ be the subalgebra of $\Uu$ spanned by elements of positive degree, i.e., non-negative along all $\alpha_i$,  and positive along at least one $\alpha_i$ ($i \in \indx_0$). 

We will also need to consider other subalgebras of $\Uu$ spanned by elements of specified degree. 
Write $\degdr_i$ for degree along the simple root $\alpha_i$. Let $Q_i$ be the sublattice of~$Q$ spanned by all the simple roots except $\alpha_i$. 
Let $\Uu_{d_i=r,+}$ be the subalgebra of $\Uu$ spanned by elements which are both: (i) of degree $r$ along $\alpha_i$, and (ii) positive $Q_i$-degree, i.e., non-negative along all the simple roots and positive along at least one simple root in $Q_i$. That is, 
\[
\Uu_{d_i=r,+} = \langle y \in \Uu \mid \degdr_i y = r;\ \forall j \neq i \in \indx_0\colon \degdr_j y \geq 0;\ \exists j \neq i\colon \degdr_j y > 0 \rangle.
\]
Moreover, set 
\begin{align}
\Uu_{d_i\geq r,+} =& \ \langle y \in \Uu \mid \degdr_i y \geq r;\ \forall j \neq i \in \indx_0\colon \degdr_j y \geq 0;\ \exists j \neq i\colon \degdr_j y > 0 \rangle, \\
\Uu_{\neq i,+} =& \ \bigoplus_{r \in \Z} \Uu_{d_i=r,+} = \langle y \in \Uu \mid \forall j \neq i \in \indx_0\colon \degdr_j y \geq 0;\ \exists j \neq i\colon \degdr_j y > 0 \rangle, \\
\Uu_+ =& \ \langle y \in \Uu \mid \forall j \in \indx_0\colon \degdr_j y \geq 0,\ \exists j \colon \degdr_j y > 0 \rangle.  
\end{align} 
We will also need the following lemma. 

\Lem \label{lem: deg of hr}
We have 
\[
\degdr e_0^{\pm} = \mp\theta. 
\]
\enlem 

\Proof
This follows directly from the explicit description of the images of $E_0$ and $F_0$ in the Drinfeld realization, see, e.g., \cite[Remark 4.7]{beck-94} and \cite[Theorem 2.2]{chari-pressley-94}. 
\enproof

\subsection{The braid group action} 
\label{subsec: br action Lusztig}

The Weyl group $W_0$ of $\g$ is generated by the simple reflections $s_i$ ($i \in \indx_0$), and acts on $P$ in the usual way. The extended affine Weyl group $\widetilde{W}$ is the semi-direct product $W_0 \ltimes P$. 
It contains the affine Weyl group $W = W_0 \ltimes Q = \langle s_i \mid i \in \indx \rangle$ as a subgroup, and can also be realized as the semi-direct product 
\eq \label{eq: ext aff weyl 2 rel}
\widetilde{W} \cong \Lambda \ltimes W,
\eneq
where $\Lambda = P/Q$ is a finite group of automorphisms of the Dynkin diagram of $\widehat{\g}$. Let $\widetilde{\mathcal{B}}$ be the braid group associated to $\widetilde{W}$. 

It is well known \cite[\S 37.1.3]{lusztig-94} that,
for each $i \in \indx$, there exists an automorphism $T_i$ of $\qla$ such that $T_i(K_\mu) = K_{s_i\mu}$ and 
\begin{alignat}{3}
T_i(E_i) =& \ -F_iK_i, \quad& T_i(E_j) =& \ \sum_{r+s = -a_{ji}} (-1)^r q_i^{-r} E_i^{(s)}E_jE_i^{(r)}, \\ 
T_i(F_i) =& \ -K_i\mi E_i, \quad& T_i(F_j) =& \ \sum_{r+s = -a_{ji}} (-1)^r q_i^{r} F_i^{(r)}F_jF_i^{(s)}, 
\end{alignat} 
for $\mu \in P$ and $i \neq j$. 
For each $\lambda \in \Lambda$, there is also an automorphism $T_\lambda$ such that $T_\lambda(E_i) = E_{\lambda(i)}$, $T_\lambda(F_i) = F_{\lambda(i)}$ and $T_\lambda(K_i) = K_{\lambda(i)}$.   
These automorphisms satisfy the relations of the braid group $\widetilde{\mathcal{B}}$. 
We extend this action to an action on $\Uu$ 
by setting 
$T_i(\KK_\mu) = \KK_{s_i\mu}$ and $T_\lambda(\KK_i) = \KK_{\lambda(i)}$. 
Given a reduced expression $w = \lambda s_{i_1} \cdots s_{i_r} \in \widetilde{W}$, one defines $T_w = T_\lambda T_{i_1} \cdots T_{i_r}$. This is independent of the choice of reduced expression.

\subsection{Quantum symmetric pairs of split affine type}

We recall the definition of quantum symmetric pair coideal subalgebras of split affine type, i.e., corresponding to Satake diagrams with no black nodes and no involution. More precisely, we consider a ``universal'' version of these algebras introduced in \cite{LW-H1, ZhangDr}, also known as universal $\imath$quantum groups. 

\Defi
Let $\Ui = \Ui(\widehat{\g})$ be the algebra generated by $\Bg_i$ and central invertible elements $\KK_i$ $(i \in \indx)$ subject to relations 
\[
-q_i \KK_i\mi\Serre_{ij}(B_i,B_j) = 
\left\{
\begin{array}{r l l}
0 & \mbox{if} & a_{ji} = 0, \\
\Bg_j  & \mbox{if} & a_{ji} = -1, \\
{[2]}_{q_i}^2[\Bg_i, \Bg_j] & \mbox{if} & a_{ji} = -2, \\ 
(1+[3]_{q_i})(B_i^2B_j+B_jB_i^2) \\
 - [4]_{q_i}(1+[2]_{q_i}^2)B_iB_jB_i + q_i^{-1}[3]_{q_i}^2\KK_iB_j & \mbox{if} & a_{ji} = -3. 
\end{array}
\right.
\]
\enDefi

It follows from \cite[Theorem 8.3]{kolb-14} that $Z(\Ui) = \C[\KK_i^{\pm1} \mid i \in \indx ]$.  
Recently, Lu--Wang \cite{lu-wang-21} and Zhang \cite{ZhangDr} found a ``Drinfeld-type" presentation of $\Ui$. 

\Defi 
Let $\DrUi$ be the $\kor$-algebra generated by $\Hg_{i,m}$ and $\Ag_{i,r}$, where $m\geq1$, $r\in\Z$, and invertible central elements $\KK_i$ ($i \in \indx_0$), $\CCC$, subject to the following relations:
\begin{align} 
[\Hg_{i,m},\Hg_{j,l}] &=0, \label{eq: grel1} \\
[\Hg_{i,m}, \Ag_{j,r}] &= \textstyle \frac{[m \cdot a_{ji}]_{q_i}}{m} (\Ag_{j,r+m}- \Ag_{j,r-m}\CCC^m), \label{eq: grel2} \\ 
[\Ag_{i,r}, \Ag_{j,s}] &= 0 \quad \mbox{if} \quad a_{ji} = 0, \label{eq: grel3} \\ 
[\Ag_{i,r}, \Ag_{j,s+1}]_{q_i^{-a_{ji}}}  -q_i^{-a_{ji}} [\Ag_{i,r+1}, \Ag_{j,s}]_{q^{a_{ji}}} &= 0  \quad \mbox{if} \quad i\neq j, \label{eq: grel4} \\ 
\label{eq: grel5}
[\Ag_{i,r}, \Ag_{i,s+1}]_{q_i^{-2}}  -q_i^{-2} [\Ag_{i,r+1}, \Ag_{i,s}]_{q_i^{2}}
&= q_i^{-2}\KK_i\CCC^r\Thg_{i,s-r+1} - q_i^{-4}\KK_i\CCC^{r+1}\Thg_{i,s-r-1}  \\ \notag
&\quad  + q_i^{-2}\KK_i\CCC^s\Thg_{i,r-s+1}  -q_i^{-4}\KK_i \CCC^{s+1}\Thg_{i,r-s-1}, 
\end{align}
and Serre relations (see \cite[(3.27)-(3.28)]{ZhangDr}), 
where $m,n\geq1$; $r,s, r_1, r_2\in \Z$ and 
\[
1+ \sum_{m=1}^\infty (q-q^{-1})\Thg_{i,m} z^m  =  \exp \left( (q-q^{-1}) \sum_{m=1}^\infty \Hg_{i,m} z^m \right).
\]
By convention, $\Thg_{i,0} = (q-q^{-1})^{-1}$ and $\Thg_{i,m} = 0$ for $m\leq -1$. 
\enDefi 

By \cite[Theorem 3.4]{ZhangDr}, there is an algebra isomorphism 
\eq
\Ui \isoto{} \DrUi.
\eneq
This isomorphism and its inverse can be described in terms of the braid group action recalled below in \S \ref{subsec: qsp br act}. For details, see \cite[(3.1)-(3.3)]{ZhangDr}. In particular, it sends
\eq
\Bg_{i} \mapsto \Ag_{i,0}, \quad \KK_i \mapsto \KK_i, \quad   \KK_\delta \mapsto \Ck, 
\eneq
for $i \in \indx_0$. 

\subsection{The QSP braid group action} \label{subsec: qsp br act}

By \cite[Lemma 2.7]{ZhangDr}, for each $i \in \indx$, there exists an automorphism $\Tbr_i$ of $\Ui$ such that $\Tbr_i(\KK_\mu) = \KK_{s_i\mu}$ and 
\[
\Tbr_i(B_j) = 
\left\{ \begin{array}{l l}
\KK_j\mi B_j & \mbox{if } i=j, \\[2pt]
B_j & \mbox{if } a_{ji} = 0, \\[2pt]
B_jB_i - q_iB_iB_j  & \mbox{if } a_{ji} = -1, \\[2pt]
{[2]_{q_i}^{-1}}(B_jB_i^2-q_i[2]_{q_i}B_iB_jB_i+q_i^2B_i^2B_j) + B_j\KK_i & \mbox{if } a_{ji} = -2, \\[2pt]
{[3]}_{q_i}\mi[2]_{q_i}\mi\big( B_jB_i^3 - q_i[3]_{q_i}B_iB_jB_i^2 +q^2[3]_{q_i}B_i^2B_jB_i \\[2pt]
-q_i^3B_i^3B_j + q_i\mi[B_j,B_i]_{q_i^3}\KK_i\big) + [B_j,B_i]_{q_i}\KK_i & \mbox{if } a_{ji} = -3, 
\end{array} \right. 
\] 
for $\mu \in \Z\indx$ and $j \in \indx$. 
For each $\lambda \in \Lambda$, there is also an automorphism $\Tbr_\lambda$ such that $\Tbr_\lambda(B_i) = B_{\lambda(i)}$ and $\Tbr_\lambda(\KK_i) = \KK_{\lambda(i)}$.   
These automorphisms define an action of the braid group $\widetilde{\mathcal{B}}$. 
We will refer to this action as the QSP braid group action to distinguish it from Lusztig's braid group action from \S \ref{subsec: br action Lusztig}. 

We will need the following lemma. 

\Lem \label{lem: br act on wt}
We have $T_w(e^{\pm}_i) = e^{\pm}_{wi}$, $\Tbr_w(B_i) = B_{wi}$, for $i \in \indx$ and $w \in \widetilde{W}$ such that $wi \in \indx$. 
\enlem

\Proof
See \cite{lusztig-94} and \cite[Lemma 2.10]{ZhangDr}. 
\enproof

\subsection{Coideal structures}

\nc{\etas}{\eta_{\mathbf{s}}}
\nc{\etacs}{\eta_{\mathbf{c}, \mathbf{s}}}
\nc{\Uic}{\mathbf{U}^\imath_{\mathbf{c}}}

By \cite[Theorem 7.1]{kolb-14}, for any $\mathbf{s} = (s_0, \cdots, s_n) \in \C^{n+1}$, 
there exists an injective algebra homomorphism 
\eq
\label{eq: Kolb emb}
\etas \colon \Ui \monoto \Uu, \quad \Bg_i \mapsto \Eg_i^- - q_i^{-2}\KK_i \Eg_i^+ \Kg_i^{-1} + \sss_i \Kg_i^{-1}, \quad  
\KK_i \mapsto \KK_i \quad 
(i \in \indx). 
\eneq 
Given $\mathbf{c} = (c_0, \cdots, c_n) \in (\C^\times)^{n+1}$, let $\Uic$ be the quotient of $\Ui$ by the two-sided ideal generated by $q_i^{-2}\KK_i - c_i$. The induced map $\etacs \colon \Uic \monoto \Uq$ 
gives $\Uic$ the structure of a right coideal subalgebra of $\Uq$, with coproduct $\Delta_{\mathbf{c},\mathbf{s}} = \Delta \circ  \eta_{\mathbf{c},\mathbf{s}}$. 
If $\mathbf{s} = (0, \cdots, 0)$, we abbreviate $\eta = \eta_{\mathbf{s}}$ and $\eta_{\mathbf{c}} = \eta_{\mathbf{c}, \mathbf{s}}$. 
Explicitly, 
\eq \label{eq: coproduct explicit} 
\Delta_{\mathbf{c},\mathbf{s}}(B_i) = 1 \otimes \eta_{\mathbf{c}}(B_i) + \etacs(B_i) \otimes K_i\mi. 
\eneq

\subsection{Rank one subalgebras} \label{subsec: rk one subalg}

\nc{\Uii}{\Ui_{[i]}} 
\nc{\Uqi}{\Uq_{[i]}} 
\nc{\Uui}{\Uu_{[i]}} 
\nc{\Usl}{U_q(L\mathfrak{sl}_2)}
\nc{\Uusl}{\widetilde{U}_q(L\mathfrak{sl}_2)}

Let us recall how the braid group actions can be used to obtain rank one subalgebras in $\Uu$ and $\Ui$. 
For $i \in \indx_0$, let $\omega'_i = \omega_i s_i$, and let $\Uqi$ be the subalgebra of $\Uq$ generated by 
\[
E_i, \ F_i, \ K_i^{\pm 1}, \ T_{\omega'_i}(E_i), \ T_{\omega'_i}(F_i), \ T_{\omega'_i}(K_i^{\pm 1}). 
\]
By \cite[Proposition 3.8]{beck-94}, there is an algebra isomorphism $\iota_i \colon \Usl \to \Uqi$ sending $q \mapsto q_i$ and 
\[
E_1 \mapsto E_i, \ F_1 \mapsto F_i, \ K_1^{\pm1} \mapsto K_i^{\pm1}, \quad 
E_0 \mapsto T_{\omega'_i}(E_i), \ F_0 \mapsto T_{\omega'_i}(F_i), \ K_0^{\pm1} \mapsto T_{\omega'_i}(K_i^{\pm1}). 
\] 
Let $\Uui$ be the subalgebra of $\Uu$ generated by $\Uqi$, $\mathbb{K}_i^{\pm1}$ and $(\KK_\delta \KK_i^{-1})^{\pm1}$. The algebra isomorphism $\iota_i$ extends to an isomorphism $\iota_i \colon \Uusl \to \Uui$ sending  
$\KK_1 \mapsto \KK_i$ and $\KK_0 \mapsto \KK_\delta \KK_i^{-1}$. 

Moreover, for $i \in \indx_0$, let $\Uii$ be the subalgebra of $\Ui$ generated by $B_i$, $\Tbr_{\omega'_i}(B_i)$, $\KK_i^{\pm1}$ and $(\KK_\delta \KK_i^{-1})^{\pm1}$. By \cite[Proposition 3.9]{lu-wang-21}, there is an algebra isomorphism $\iota \colon \Ui(\widehat{\mathfrak{sl}}_2) \to \Uii$ sending $q \mapsto q_i$ and 
\eq \label{eq: rank 1 q onsager}
B_1 \mapsto B_i, \quad B_0 \mapsto \Tbr_{\omega'_i}(B_i), \quad \KK_1 \mapsto \KK_i, \quad \KK_0 \mapsto \KK_\delta \KK_i^{-1}. 
\eneq 
We also note that $A_{1,-1}$ is sent to $\Tbr_{\omega_i}(B_i)$.

\subsection{Factorization and coproduct in rank one}

\nc{\icar}{\mathcal{H}^\imath}
\nc{\Tser}{\pmb{\grave{\Theta}}}
\nc{\xis}{\xi_{\mathbf{s}}}

Let $\car$ denote the commutative subalgebra of $\Uu$ generated by the coefficients $\phi^{\pm}_{i,\pm r}$ of the series $\pmb{\phi}^{\pm}_i(z)$. We will refer to these coefficients as the Drinfeld--Cartan operators, and to $\car$ as the Drinfeld--Cartan subalgebra of $\Uu$. 
Similarly, we refer to the coefficients of the series 
\[
\pmb{\grave{\Theta}}_i(z) = (q_i-q_i\mi)\frac{1-q_i^{-2}\CCC z^2}{1-\CCC z^2} \sum_{r \geq 0} \Theta_{i,r} z^r
\]
as $\imath$Drinfeld--Cartan operators\footnote{One has $\pmb{\grave{\Theta}}_i(z) = (q_i-q_i\mi) \acute{\Theta}_i(z)$, where $\acute{\Theta}_i(z) = \sum_{r\geq 0} \acute{\Theta}_{i,r} z^r$ is the series from \cite[(2.10)]{lu-wang-21}, whose coefficients are the imaginary root vectors from \cite{bas-kol-20}, rescaled by $-q_i^{-2}$.}, 
and to the commutative subalgebra $\icar$ generated by them as the $\imath$Drinfeld--Cartan subalgebra of $\Ui$. It is natural to ask how the two Drinfeld--Cartan subalgebras are related under the monomorphism \eqref{eq: Kolb emb}. This question was answered in \cite{Przez-23} for quantum symmetric pairs of split affine type $\mathsf{A}$. Here we will need the rank one case of that result, i.e., when $\Ui$ is isomorphic to the (universal) $q$-Onsager algebra. 

\Thm[{\cite[Theorem 7.5]{Przez-23}}] \label{thm: rank 1 factorization}
Let $\widehat{\g} = \widehat{\mathfrak{sl}}_2$. 
The series $\Tser(z)$ admits the following factorization 
\begin{align} 
\etas(\Tser(z)) \equiv& \ \xis(\Tser(z)) \cdot \pmb{\phi}^-(z\mi)\pmb{\phi}^+(\CCC z) \quad  \mod \fext{\Uu_+}{z}, \\ 
\Delta_{\mathbf{s}}(\Thgsr(z)) \equiv& \ \eta_{\mathbf{s}}(\Thgsr(z)) \otimes \eta(\Thgsr(z)) \quad  \mod \fext{\Uu  \otimes \Uu_{+}}{z}, 
\end{align} 
where $\xis = \varepsilon \circ \etas$. 
\enthm


\section{Reduced expressions for fundamental weights}
\label{sec-comb}

The goal of this section is to obtain explicit reduced expressions for the extended affine Weyl group elements corresponding to the fundamental weights in types $\mathsf{A}_n^{(1)}, \mathsf{B}_n^{(1)}$, $\mathsf{C}_n^{(1)}, \mathsf{D}_n^{(1)}$.

First, we need to introduce some notations. 
Let $w_0$ be the longest element of the finite Weyl group $W_0$ of $\g$.
Given $i \in \indx_0$, let $w_i$ be the longest element of the Weyl group corresponding to the Dynkin diagram of $\g$ with the node labelled by $i$ removed. 
Given $1 \leq k < l$, let 
\[ [k,l] = s_k s_{k+1} \cdots s_{l}, \quad [l,k] = [k,l]^{-1}.\]
Given the highest root $\theta = \sum_{i \in \indx_0} c_i \alpha_i$, let $J = \{ i \in \indx_0 \mid c_i = 1\}$. 

We will use the following two results. 

\Pro[{\cite[Ch.\ VI \S 2.3, Proposition 6]{bourbaki}, \cite[Proposition 1.18]{IM65}}] \label{pro: bourbaki}
There is a bijection
$
J  \longrightarrow  \Lambda {-} \{1\}
$
sending $i \mapsto \omega_i w_i w_0$. 
\enpro 


\Pro[{\cite[Proposition 1.23]{IM65}}] \label{pro: length wt}
The length of the fundamental weights is given by
\[
\ell(\omega_i) = \sum_{\beta \in \Delta^+}(\beta, \omega_i),
\]
where $\Delta^+$ is the set of the positive roots of $\g$. 
\enpro 


\begin{table}
\centering
\begin{tabular}{|l|c|}
\hline
\begin{tabular}{l}
$\mathsf{A}_{n}^{(1)}$ 
\end{tabular}
 & 
 \begin{tabular}{l}
 \begin{tikzpicture}[scale=0.5]
\draw (0 cm,0) -- (2 cm,0);
\draw (2 cm,0) -- (4 cm,0);
\draw[dashed] (4 cm,0) -- (6 cm,0);
\draw (6 cm,0) -- (8 cm,0);
\draw (8 cm,0) -- (10 cm,0);
\draw (0 cm,0) -- (5.0 cm, 1.2 cm);
\draw (5.0 cm, 1.2 cm) -- (10 cm, 0);
\draw[fill=white] (0 cm, 0 cm) circle (.25cm) node[below=4pt]{$1$};
\draw[fill=white] (2 cm, 0 cm) circle (.25cm) node[below=4pt]{$2$};
\draw[fill=white] (4 cm, 0 cm) circle (.25cm) node[below=4pt]{};
\draw[fill=white] (6 cm, 0 cm) circle (.25cm) node[below=4pt]{};
\draw[fill=white] (8 cm, 0 cm) circle (.25cm) node[below=4pt]{$n-1$};
\draw[fill=white] (10 cm, 0 cm) circle (.25cm) node[below=4pt]{$n$};
\draw[fill=white] (5.0 cm, 1.2 cm) circle (.25cm) node[anchor=south east]{$0$};
\end{tikzpicture} \end{tabular} \\
\hline
\begin{tabular}{l}
$\mathsf{B}_{2}^{(1)}$ 
\end{tabular}
 &
 \begin{tabular}{l} \\
\begin{tikzpicture}[scale=0.5]
\draw (0, 0.1 cm) -- +(2 cm,0);
\draw (0, -0.1 cm) -- +(2 cm,0);
\draw[shift={(1.2, 0)}, rotate=0] (135 : 0.45cm) -- (0,0) -- (-135 : 0.45cm);
{
\pgftransformxshift{2 cm}
\draw (0 cm,0) -- (0 cm,0);
\draw (0 cm, 0.1 cm) -- +(2 cm,0);
\draw (0 cm, -0.1 cm) -- +(2 cm,0);
\draw[shift={(0.8, 0)}, rotate=180] (135 : 0.45cm) -- (0,0) -- (-135 : 0.45cm);
\draw[fill=white] (0 cm, 0 cm) circle (.25cm) node[below=4pt]{$2$};
\draw[fill=white] (2 cm, 0 cm) circle (.25cm) node[below=4pt]{$1$};
}
\draw[fill=white] (0 cm, 0 cm) circle (.25cm) node[below=4pt]{$0$};
\end{tikzpicture} \end{tabular} \\
\hline
\begin{tabular}{l}
$\mathsf{B}_{n}^{(1)}$, $n \ge 3$
\end{tabular}
  &
\begin{tabular}{l}
\begin{tikzpicture}[scale=0.5]
\draw (0,0.7 cm) -- (2 cm,0);
\draw (0,-0.7 cm) -- (2 cm,0);
\draw (2 cm,0) -- (4 cm,0);
\draw[dashed] (4 cm,0) -- (6 cm,0);
\draw (6 cm,0) -- (8 cm,0);
\draw (8 cm, 0.1 cm) -- +(2 cm,0);
\draw (8 cm, -0.1 cm) -- +(2 cm,0);
\draw[shift={(9.2, 0)}, rotate=0] (135 : 0.45cm) -- (0,0) -- (-135 : 0.45cm);
\draw[fill=white] (0 cm, 0.7 cm) circle (.25cm) node[left=3pt]{$0$};
\draw[fill=white] (0 cm, -0.7 cm) circle (.25cm) node[left=3pt]{$1$};
\draw[fill=white] (2 cm, 0 cm) circle (.25cm) node[below=4pt]{$2$};
\draw[fill=white] (4 cm, 0 cm) circle (.25cm) node[below=4pt]{};
\draw[fill=white] (6 cm, 0 cm) circle (.25cm) node[below=4pt]{};
\draw[fill=white] (8 cm, 0 cm) circle (.25cm) node[below=4pt]{$n-1$};
\draw[fill=white] (10 cm, 0 cm) circle (.25cm) node[below=4pt]{$n$};
\end{tikzpicture} \end{tabular} \\
\hline
\begin{tabular}{l}
$\mathsf{C}_{n}^{(1)}$, $n \ge 3$
\end{tabular}
  &
\begin{tabular}{l} \\
\begin{tikzpicture}[scale=0.5]
\draw (0, 0.1 cm) -- +(2 cm,0);
\draw (0, -0.1 cm) -- +(2 cm,0);
\draw[shift={(1.2, 0)}, rotate=0] (135 : 0.45cm) -- (0,0) -- (-135 : 0.45cm);
{
\pgftransformxshift{2 cm}
\draw (0 cm,0) -- (2 cm,0);
\draw (2 cm,0) -- (4 cm,0);
\draw[dashed] (4 cm,0) -- (6 cm,0);
\draw (6 cm,0) -- (8 cm,0);
\draw (8 cm, 0.1 cm) -- +(2 cm,0);
\draw (8 cm, -0.1 cm) -- +(2 cm,0);
\draw[shift={(8.8, 0)}, rotate=180] (135 : 0.45cm) -- (0,0) -- (-135 : 0.45cm);
\draw[fill=white] (0 cm, 0 cm) circle (.25cm) node[below=4pt]{$1$};
\draw[fill=white] (2 cm, 0 cm) circle (.25cm) node[below=4pt]{$2$};
\draw[fill=white] (4 cm, 0 cm) circle (.25cm) node[below=4pt]{};
\draw[fill=white] (6 cm, 0 cm) circle (.25cm) node[below=4pt]{};
\draw[fill=white] (8 cm, 0 cm) circle (.25cm) node[below=4pt]{$n-1$};
\draw[fill=white] (10 cm, 0 cm) circle (.25cm) node[below=4pt]{$n$};
}
\draw[fill=white] (0 cm, 0 cm) circle (.25cm) node[below=4pt]{$0$};
\end{tikzpicture} \end{tabular} \\
\hline
\begin{tabular}{l}
$\mathsf{D}_{n}^{(1)}$, $n \ge 4$
\end{tabular}
  & 
\begin{tabular}{l}
\begin{tikzpicture}[scale=0.5]
\draw (0,0.7 cm) -- (2 cm,0);
\draw (0,-0.7 cm) -- (2 cm,0);
\draw (2 cm,0) -- (4 cm,0);
\draw[dashed] (4 cm,0) -- (6 cm,0);
\draw (6 cm,0) -- (8 cm,0);
\draw (8 cm,0) -- (10 cm,0.7 cm);
\draw (8 cm,0) -- (10 cm,-0.7 cm);
\draw[fill=white] (0 cm, 0.7 cm) circle (.25cm) node[left=3pt]{$0$};
\draw[fill=white] (0 cm, -0.7 cm) circle (.25cm) node[left=3pt]{$1$};
\draw[fill=white] (2 cm, 0 cm) circle (.25cm) node[below=4pt]{$2$};
\draw[fill=white] (4 cm, 0 cm) circle (.25cm) node[below=4pt]{};
\draw[fill=white] (6 cm, 0 cm) circle (.25cm) node[below=4pt]{};
\draw[fill=white] (8 cm, 0 cm) circle (.25cm) node[below=4pt]{$n-2$};
\draw[fill=white] (10 cm, 0.7 cm) circle (.25cm) node[right=3pt]{$n$};
\draw[fill=white] (10 cm, -0.7 cm) circle (.25cm) node[right=3pt]{$n-1$};
\end{tikzpicture} \end{tabular}
\\
\hline
\end{tabular}
\caption{Affine Dynkin diagrams of types $\mathsf{A}_n^{(1)}$, $\mathsf{B}_n^{(1)}$, $\mathsf{C}_n^{(1)}$, $\mathsf{D}_n^{(1)}$.}
\label{table:Dynkin diagrams of type ABCD}
\end{table}

\subsection{Type $\mathsf{A}_n^{(1)}$} 

The highest root is $\theta = \alpha_1 + \alpha_2 + \cdots + \alpha_n$ and $J = \indx_0$.
The fundamental group $\Lambda$ is cyclic of order $n+1$. It is generated by, e.g., the automorphism $\pi$ of the affine Dynkin diagram sending $i \mapsto i+1 \mod n+1$.

\Pro \label{pro: type A wt}
Let $\widetilde{W}$ be the extended affine Weyl group of type $A_n^{(1)}$, for $n \geq 1$. Then 
\[
\omega_i = \pi^i [n - i+1 ,n] \cdots[2,i+1][1,i]. 
\] 
This is a reduced expression. 
\enpro

\Proof
See \cite[\S 4.5]{Lus83} or \cite[Proposition 9.2]{Przez-23}. 
\enproof

\subsection{Type $\mathsf{D}_n^{(1)}$} 

The highest root is $\theta = \alpha_1 + \alpha_{n-1} + \alpha_n + 2(\alpha_2 + \cdots + \alpha_{n-2})$ and $J = \{1, n{-}1, n\}$. The fundamental group $\Lambda = \{ 1, \pi_1, \pi_{n-1}, \pi_n \}$ depends on the parity of $n$. If $n$ is even then $\Lambda \cong \Z/2\Z \times \Z/2\Z$, and if $n$ is odd then $\Lambda \cong \Z/4\Z$. Explicit descriptions of the corresponding diagram automorphisms can be found in \cite[Ch.\ VI \S 4.8]{bourbaki}. 
Given $1 \leq m \leq n-2$, set 
\[
r_m = s_n [m,n{\shortminus}1][m{\shortminus}1,n{\shortminus}2].
\]

\Pro \label{pro: type D funda weights}
Let $\widetilde{W}$ be the extended affine Weyl group of type $\mathsf{D}_n^{(1)}$, for $n \geq 4$. If $i$ is even, then 
\[
\omega_i = \left\{
\begin{array}{r l }
\big(s_0 [2,n{\shortminus}1][1,n{\shortminus}2]s_n\big)^{\frac{i}{2}} [{n{\shortminus}i},{n{\shortminus}1}] \cdots [1,i] & \quad \mbox{if }\  2 \leq i \leq n{\shortminus}2, \\ 
\pi_{n{\shortminus}1} r_{n{\shortminus}2}r_{n{\shortminus}4} \cdots r_3 s_n [1,{n{\shortminus}1}] 
& \quad \mbox{if }\  i = n{\shortminus}1, \\  
\pi_{n} r_{n{\shortminus}2} r_{n{\shortminus}4} \cdots r_2 s_n 
& \quad \mbox{if }\  i = n. \\  
\end{array}
\right.
\]
If $i$ is odd, then 
\[
\omega_i = \left\{
\begin{array}{r l }
\pi_1[1,n{\shortminus}2]s_n\big(s_0[2,n{\shortminus}1][1,n{\shortminus}2]s_n\big)^{\frac{i{\shortminus}1}{2}}[{n{\shortminus}i},{n{\shortminus}1}] \cdots [1,i] & \quad \mbox{if }\  1 \leq i \leq n{\shortminus}2, \\ 
\pi_{n{\shortminus}1} s_{n{\shortminus}1}s_{n{\shortminus}2} r_{n{\shortminus}3} r_{n{\shortminus}5} \cdots r_3 s_n [1,{n{\shortminus}1}] 
& \quad \mbox{if }\  i = n{\shortminus}1, \\  
\pi_{n} s_{n{\shortminus}1}s_{n{\shortminus}2}r_{n{\shortminus}3} r_{n{\shortminus}5} \cdots r_2 s_n 
& \quad \mbox{if }\  i = n. \\  
\end{array}
\right. 
\]
Moreover, the formulae above yield reduced expressions.
\enpro 

\Proof 
If $i \in J = \{1,n{-}1,n\}$, one may use Proposition \ref{pro: bourbaki} and argue in a similar way as in Proposition \ref{pro: type A wt}. 
For more general $i \in \indx_0$, consider $\widetilde{W}$ as the group of affine transformations of $\R^n$. 
We show that elements on the RHS of the equations above act in the same way as $\omega_i$. 
Indeed, the element $[{n{\shortminus}i},{n{\shortminus}1}] \cdots [1,i]$ acts as a cyclic permutation of the coordinates (shifting the indices by $n-i$). It is also easy to explicitly verify that $\big(s_0 [2,n{\shortminus}1][1,n{\shortminus}2]s_n\big)^{\frac{i}{2}}$ (in the even case) and $\pi_1\big([1,n{\shortminus}2]s_ns_0[2,n{\shortminus}1]\big)^{\frac{i{\shortminus}1}{2}}[1,n{\shortminus}2]s_n$ (in the odd case) act as cyclic permutations of the coordinates (shifting the indices by $i$) followed by an affine shift by $\omega_i$.  

We will prove the reducedness of our expressions for $\omega_i$ by induction on $n$. The base case of $\mathsf{D}_4$ can be checked directly. Now assume that $n \ge 5$, and that the result is true for $\mathsf{D}_{n-1}$. The set of positive roots of type $\mathsf{D}_{n}$ ($n \ge 5$) is the union of the set of positive roots of type $\mathsf{D}_{n-1}$ (with the indices shifted by $1$), and the set consisting of 
the roots 
\begin{align} \label{eq:second part of positive roots in Dn}
\sum_{j=1}^l \alpha_j \ (1 \leq l \leq n), \quad \alpha_n+\sum_{j=1}^{n-2}\alpha_j, \quad \sum_{j=1}^{n} \alpha_j + \sum_{j=2}^{l} \alpha_{n-j} \ (1 \leq l \leq n-2).
\end{align}

Let $1 \le i \le n-2$. We will show that $\ell(\omega_i) = i(2n-i-1)$. By Proposition \ref{pro: length wt}, it suffices to check that the sum of the coefficients of $\alpha_i$ in positive roots of type $\mathsf{D}_n$ (written in terms of the basis of simple roots) is equal to $i(2n-i-1)$. This sum is equal to the sum of the coefficients of $\alpha_{i-1}$ in the positive roots of type $\mathsf{D}_{n-1}$ and the sum of the coefficients of $\alpha_i$ in the roots in \eqref{eq:second part of positive roots in Dn}. By induction, the former is equal to $(i-1)(2n-2-i)$. Therefore, the total sum is $ (i-1)(2n-2-i) + (n-i+1) + 1 + (n-3) + (i-1) = i(2n-i-1)$. 

Next, let $i=n$. We will check that $\ell(\omega_n)=n(n-1)/2$. By induction, the sum of the coefficients of $\alpha_{n-1}$ in the positive roots of type $\mathsf{D}_{n-1}$ is $(n-1)(n-2)/2$. Therefore the total sum of the coefficients of $\alpha_{n}$ in the positive roots of type $\mathsf{D}_{n}$ is $(n-1)(n-2)/2+n-1 = n(n-1)/2$.

Finally, let $i=n-1$. We will check that $\ell(\omega_n)=n(n-1)/2$. By induction, the sum of the coefficients of $\alpha_{n-2}$ in the positive roots of type $\mathsf{D}_{n-1}$ is $(n-1)(n-2)/2$. Therefore the total sum of the coefficients of $\alpha_{n-1}$ in the positive roots of type $\mathsf{D}_{n}$ is $(n-1)(n-2)/2+n-1 = n(n-1)/2$. 

In each of the three cases, the length of $\omega_i$ equals the number of letters in the expressions for $\omega_i$ we found, implying that these expressions are reduced. 
\enproof

\subsection{Type $\mathsf{B}_n^{(1)}$} 
The highest root is $\theta=\alpha_1+ 2\sum_{i=2}^n \alpha_i$ and $J=\{1\}$. The fundamental group is $\Lambda = \{ 1, \pi_1\}$, where $\pi_1$ interchanges $\alpha_0$ and $\alpha_1$, and fixes the other $\alpha_j$'s.

\Pro
Let $\widetilde{W}$ be the extended affine Weyl group of type $\mathsf{B}_n^{(1)}$, for $n \geq 2$. If $i$ is odd, then\footnote{If $i=n$, the expression below reads $\pi_1 [1,n] \big(s_0[2,n][1,n]\big)^{\frac{n-1}{2}}$.}
\begin{align*}
\omega_i = \pi_1 [1,n] \big(s_0[2,n][1,n]\big)^{\frac{i-1}{2}} [{n{\shortminus}i},{n{\shortminus}1}] \cdots [1,i].
\end{align*}
If $i$ is even, then\footnote{If $i=n$, the expression below reads $\big(s_0[2,n][1,n]\big)^{\frac{n}{2}}$.}
\begin{align*} 
\omega_i = \big(s_0[2,n][1,n]\big)^{\frac{i}{2}}  [{n{\shortminus}i},{n{\shortminus}1}] \cdots [1,i].
\end{align*}
Moreover, the formulae above yield reduced expressions.
\enpro  

\begin{proof}
Given $I \subset \{1, n\}$, let $w_0^I$ denote the longest word in the parabolic subgroup $W_I$. 

First consider the case of $i=1$. By Proposition \ref{pro: bourbaki}, $\omega_1 = \pi_1 w_0 w_1$. We have $w_0 = w_0^{[1,n-1]} (s_n (s_{n-1}s_n) \cdots [2,n])^{-1}$ and $w_1 = s_n (s_{n-1}s_n) \cdots [2,n] w_0^{[2,n-1]}$. Therefore $\pi_1 w_0 w_1 = \pi_1 w_0^{[1,n-1]} w_0^{[2,n-1]} = \pi_1 [1,n] [n-1,1]$. 

For $i \neq 1$, we use the same strategy as in the proof of Proposition \ref{pro: type D funda weights}, i.e., show that the expressions we found define the same affine transformations of $\R^n$. More precisely, we need to show that each $\omega_i$ (as an element of the extended affine Weyl group of type $\mathsf{B}_n$) acts as a translation by the $i$-th fundamental weight of type $\mathsf{C}_n$. 
Let us first recall some facts about the root systems of types $\mathsf{B}_n$ and $\mathsf{C}_n$. 
The simple roots in the root system of type $\mathsf{B}_n$ are
\begin{align*}
\alpha_j = e_j - e_{j+1}\ (1 \leq j \leq n-1), \ \ \alpha_n=e_n. 
\end{align*}
We have $\theta=\alpha_1+2\sum_{i=2}^n \alpha_i = e_1+e_2=\theta^{\vee}$,   
\begin{alignat}{3}
s_j(e_j) =& \ e_{j+1} \ && (1 \leq j \leq n-1), \quad & s_n(e_n) &= \ -e_n, \\
s_j(e_{j+1}) =& \ e_j \ && (1 \leq j \leq n-1), \quad & s_0(e_j) &= \ t_{\theta^\vee} r_\theta(e_j), 
\end{alignat}
where $t_{\theta^\vee}$ is the translation by ${\theta^\vee}$; $r_\theta$ is the reflection by $\theta$; 
and in all other cases $s_j(e_{j'}) = e_{j'}$. 
Finally, recall that the fundamental weights of type $\mathsf{C}_n$ are $\omega_i^* = \sum_{j=1}^i e_j$ ($1 \leq i \leq n$). 

Observe that $[n-i,n-1] \cdots [2,i+1] [1,i]$ acts by cyclically shifting the coordinates by $n-i$. 
Let $i$ be even. We have  
\begin{align*}
& (s_0[2,n][1,n])^{\frac{i}{2}} [n-i, n-1] \cdots [2,i+1] [1,i] (e_j) = (s_0[2,n][1,n])^{\frac{i}{2}} (e_{n-i+j}) \\
& =  (s_0[2,n][1,n])^{\frac{i}{2}-1} t_{\theta} (e_{j-i+2}) =  (s_0[2,n][1,n])^{\frac{i}{2}-1} (e_{j-i+2} + e_1 + e_2) \\
& =  (s_0[2,n][1,n])^{\frac{i}{2}-2} (e_{j-i+4} + e_1+e_2 + e_3 + e_4) = e_j + \sum_{r=1}^i e_r = e_j + \omega_i^*.
\end{align*}
 
Let $i$ be odd. We have 
\begin{align*}
& \pi_1 ([1,n]s_0[2,n])^{\frac{i-1}{2}}[1,n] [n-i,n-1] \cdots [2,i+1][1,i] (e_j) \\
& = \pi_1  ([1,n]s_0[2,n])^{\frac{i-1}{2}} [1,n] (e_{n-i+j}) = \pi_1 [1,n] (s_0[2,n][1,n])^{\frac{i-1}{2}} (e_{n-i+j}) \\
& = \pi_1 [1,n] (s_0[2,n][1,n])^{\frac{i-3}{2}} (e_{j-i+2} + e_1+e_2) = \pi_1 [1,n] (e_{j-1} + \sum_{r=1}^{i-1} e_r) \\
& = t_{\omega_1^*}  [1,n-1] (e_{j-1} + \sum_{r=1}^{i-1} e_r) = t_{\omega_1^*} (e_{j} + \sum_{r=2}^i e_r ) = e_{j} + \omega_i^*,
\end{align*} 
where in the passage from the third to the fourth lines we used the already established formula for $\omega_1$. 

We will now prove, by induction on $n$, that the expressions for $\omega_i$ we produced are reduced. The base case $n=2$ can be checked directly. So assume that $n \ge 3$ and the result is true for $\mathsf{B}_{n-1}$. The set of positive roots of type $\mathsf{B}_{n}$ ($n \ge 3$) is the union of the set of positive roots of type $\mathsf{B}_{n-1}$ (with indices shifted by $1$) and the set consisting of 
the roots 
\begin{align} \label{eq:second part of positive roots in Bn}
\sum_{j=1}^l \alpha_j \ (1 \leq l \leq n), \quad \sum_{j=1}^{n} \alpha_j + \sum_{j=1}^{l} \alpha_{n-j+1} \ (1 \leq l \leq n-1).
\end{align} 

Let $1 \le i \le n$. We will check that $\ell(\omega_i) = i(2n-i)$. By Proposition \ref{pro: length wt}, it suffices to check that the sum of the coefficients of $\alpha_i$ in positive roots of type $\mathsf{B}_n$ is equal to $i(2n-i)$. This sum is equal to the sum of the coefficients of $\alpha_{i-1}$ in positive roots of type $\mathsf{B}_{n-1}$ and the sum of the coefficients of $\alpha_i$ in the positive roots in \eqref{eq:second part of positive roots in Bn}. By induction, the former is equal to $(i-1)(2n-2-i+1)$. Therefore, the total sum is $ (i-1)(2n-2-i+1) + (n-i+1) + (n-1) + (i-1) = i(2n-i)$. This is also the length of the expressions for $\omega_i$ we found, implying they are indeed reduced. 
\end{proof}

\subsection{Type $\mathsf{C}_n^{(1)}$}

The highest root is $\theta=\alpha_n + 2\sum_{i=1}^{n-1} \alpha_i$ and $J=\{n\}$. The fundamental group is $\Lambda = \{ 1, \pi_n\}$, where $\pi_n$ interchanges $\alpha_k$ and $\alpha_{n-k}$. 

\Pro
Let $\widetilde{W}$ be the extended affine Weyl group of type $\mathsf{C}_n^{(1)}$, for $n \geq 2$. If $i=n$, then 
\begin{align} \label{eq: om n C}
\omega_n = \pi_n s_n[n-1,n] \cdots [1,n].
\end{align}
If $i<n$, then 
\begin{align*}
\omega_i = (s_0 [1,n] )^i [n-i,n-1] \cdots [2,i+1][1,i]. 
\end{align*}
Moreover, the formulae above yield reduced expressions.
\enpro 

\begin{proof}
We follow the same strategy as in the proofs for types $\mathsf{B}_n^{(1)}$ and $\mathsf{D}_n^{(1)}$. 
The simple roots in the root system of type $\mathsf{C}_n$ are
\begin{align*} 
\alpha_j = e_j - e_{j+1} \ (1 \leq j \leq n-1), \ \ \alpha_n = 2 e_n.
\end{align*}
The highest root is $\theta = 2 e_1 = e_n + 2 \sum_{i=1}^{n-1} e_i$ and $\theta^{\vee} = e_1$. We have  
\begin{alignat}{3}
 s_j(e_j) =& \ e_{j+1} \ &&(1 \leq j \leq n-1), \quad  & s_n(e_n) =& \ -e_n, \\
 s_j(e_{j+1}) =& \ e_j  \ &&(1 \leq j \leq n-1), \quad & s_0(e_j) =& \ t_{\theta^\vee} r_\theta(e_j), 
\end{alignat}
and in all other cases $s_j(e_{j'}) = e_{j'}$. 
Moreover, the fundamental weights of the dual root system $\mathsf{B}_n$ are $\omega_i^* = \sum_{j=1}^i e_j$ ($i<n$), and $\omega_n^* =\frac{1}{2} \sum_{j=1}^n e_j$. 

The case of $i=n$ can simply be handled using Proposition \ref{pro: bourbaki}, i.e., $\omega_n = \pi_n w_0 w_n = \pi_n s_n(s_{n-1}s_n) \cdots [1,n]$. 
Next, suppose that $i<n$. We have 
\begin{align*}
& (s_0 [1,n])^i [n-i,n-1] \cdots [2,i+1][1,i] (e_j) = (s_0 [1,n])^i (e_{n-i+j}). 
\end{align*}
If $j=i$, then 
\begin{align*}
& (s_0 [1,n])^i (e_{n}) = (s_0 [1,n])^{i-1} t_{\theta^{\vee}}r_{\theta} (-e_1) = (s_0 [1,n])^{i-1} (2e_1) \\
& = (s_0 [1,n])^{i-2} (e_1+2e_2) = (s_0 [1,n])^{i-3} (e_1+e_2+2e_3) = e_i + \sum_{r=1}^i  e_r = \omega_i^* + e_i. 
\end{align*}
If $i<j$, then $j-i+1 \ge 2$ and
\begin{align*}
& (s_0 [1,n])^i (e_{n-i+j}) = (s_0 [1,n])^i (e_{j-i}) = (s_0 [1,n])^{i-1} t_{\theta^{\vee}} r_{\theta} (e_{j-i+1}) \\
& = (s_0 [1,n])^{i-1} (e_1+e_{j-i+1}) = (s_0 [1,n])^{i-2} (e_1+e_2+e_{j-i+2}) = e_j + \sum_{r=1}^i = e_j + \omega_i^*. 
\end{align*}
If $j<i$, then 
\begin{align*}
&  (s_0 [1,n] )^i (e_{n-i+j}) =  (s_0 [1,n] )^{i-1} (e_1+e_{n-i+j+1})  = (s_0 [1,n] )^{i-2} (e_1+e_2+e_{n-i+j+2}) \\
& = (s_0 [1,n] )^{j} ( e_n + \sum_{r=1}^{i-j} e_r ) = (s_0 [1,n] )^{j-1}  (e_1 + \sum_{r=1}^{i-j+1} e_r ) =e_j+ \sum_{r=1}^{i} e_r =e_j + \omega_i^*. 
\end{align*}

We will now prove, by induction on $n$, that the expressions for $\omega_i$ we found are reduced. The base case of $\mathsf{C}_2$ can be checked directly. So assume that $n \ge 3$ and the result is true for $\mathsf{C}_{n-1}$. The set of positive roots of type $\mathsf{C}_{n}$ ($n \ge 3$) is the union of the set of positive roots of type $\mathsf{C}_{n-1}$ (with the indices by $1$) and the set consisting of 
the roots 
\begin{align} \label{eq:second part of positive roots in Cn}
\sum_{j=1}^l \alpha_j \ (1 \leq l \leq n), \quad \sum_{j=1}^{n} \alpha_j + \sum_{j=1}^{l} \alpha_{n-j} \ (1 \leq l \leq n-1).
\end{align} 
We will check that $\ell(\omega_i) = i(n+1)$ if $i<n$ and $\ell(\omega_n)=n(n+1)/2$. By Proposition \ref{pro: length wt}, it suffices to check that the sum of coefficients of $\alpha_i$ in positive roots of type $\mathsf{C}_n$ is equal to $i(n-1)$ if $i<n$, and to $n(n+1)/2$ if $i=n$. This sum is equal to the sum of the coefficients of $\alpha_{i-1}$ in positive roots of type $\mathsf{C}_{n-1}$ and the sum of the coefficients of $\alpha_i$ in the positive roots in (\ref{eq:second part of positive roots in Cn}). 

Let $i=n$. By induction hypothesis, the sum of the coefficients of $\alpha_{n-1}$ in positive roots of type $\mathsf{C}_{n-1}$ is equal to $(n-1)n/2$. Therefore, the total sum is $ (n-1)n/2  + n = n(n+1)/2$.  
Now let $i<n$. By induction hypothesis, the sum of the coefficients of $\alpha_{i-1}$ in positive roots of type $\mathsf{C}_{n-1}$ is equal to $(i-1)n$. Therefore, the total sum is $ (i-1)n + (n-i+1) + (n-1) + i = i(n+1)$, as claimed. 
This is also the length of the expressions for $\omega_i$ we found, implying they are indeed reduced. 
\end{proof}

\section{Auxiliary calculations - root combinatorics} 
\label{sec: root combin}

We will need the following lemmas to compute the braid group action in \S \ref{sec: weak com}, with the aid of Lemma \ref{lem: br act on wt}. Their proofs are based on straightforward calculations, so we omit most of the details. 

\subsection{Type $\mathsf{D}_n^{(1)}$} 
We consider the cases of weights $\omega_1, \cdots, \omega_{n-2}$ and $\omega_{n-1}, \omega_n$ separately. 

\subsubsection{Weights $\omega_1, \cdots, \omega_{n-2}$}

\Lem \label{lem: aux root 1}
We have
\begin{align}
s_0 [2,n{\shortminus}1][1,n{\shortminus}2]s_n \cdot \alpha_j =& \ 
\left\{
\begin{array}{r l }
\alpha_0 + \alpha_1 + 2\alpha_2 + \alpha_3 & \quad \mbox{if }\ j=0, \\
\alpha_{j+2} & \quad \mbox{if }\ 1 \leq j \leq n-3, \\
\alpha_1 & \quad \mbox{if }\ j = n-1, \\
\alpha_0 + \alpha_2 + \cdots + \alpha_{n-2} + \alpha_n & \quad \mbox{if }\ j = n-2, \\
-\theta - 2\alpha_0 & \quad \mbox{if }\ j = n, 
\end{array}
\right. \\ 
\pi_1[1,n{\shortminus}2]s_n \cdot \alpha_j  =& \ 
\left\{
\begin{array}{r l }
\alpha_0 + \alpha_1 + \alpha_2 & \quad \mbox{if }\ j=0, \\
\alpha_{j+1} & \quad \mbox{if }\ 1 \leq j \leq n-3, \\
\alpha_0 + \alpha_2 + \cdots + \alpha_{n-2} + \alpha_{n-1} & \quad \mbox{if }\ j = n-1, \\
\alpha_n & \quad \mbox{if }\ j = n-2, \\
-(\alpha_0 + \alpha_2 + \cdots + \alpha_{n-2} + \alpha_n) & \quad \mbox{if }\ j = n, 
\end{array}
\right. 
\end{align}
\begin{align}
(s_0 [2,n{\shortminus}1][1,n{\shortminus}2]s_n)^2 \cdot \alpha_{n-2} =& \ \alpha_2, \\
\pi_1[1,n{\shortminus}2]s_n(s_0 [2,n{\shortminus}1][1,n{\shortminus}2]s_n) \cdot \alpha_{n-2} =& \ \alpha_1. 
\end{align}
\enlem

\Proof
The lemma follows by direct calculation. 
\enproof

Set
\begin{align}
\tilde{\alpha}_k =& \ \left\{
\begin{array}{r l}
\alpha_k & \quad \mbox{if } \ 1 \leq k \leq n-1, \\
\alpha_0 + \alpha_2 + \cdots + \alpha_{n-2} + \alpha_n & \quad \mbox{if }\ k = 0. 
\end{array}
\right. \\
\hat{\alpha}_k =& \ \left\{
\begin{array}{r l}
\alpha_k & \quad \mbox{if } \ 1 \leq k \leq n-2, \\ 
\alpha_{n} & \quad \mbox{if } \ k = n-1, \\
\alpha_0 + \alpha_2 + \cdots + \alpha_{n-2} + \alpha_{n-1} & \quad \mbox{if }\ k = 0. 
\end{array}
\right.
\end{align}
For $1 \leq i \leq i-2$, set 
\[
\zeta_i = \left\{
\begin{array}{r l}
\big(s_0 [2,n{\shortminus}1][1,n{\shortminus}2]s_n\big)^{\frac{i}{2}} & \quad \mbox{if} \ i \ \mbox{is even}, \\
\pi_1[1,n{\shortminus}2]s_n\big(s_0[2,n{\shortminus}1][1,n{\shortminus}2]s_n\big)^{\frac{i{\shortminus}1}{2}} 
& \quad \mbox{if} \ i \ \mbox{is odd}. 
\end{array}
\right.
\] 


\Lem \label{lem: action on sroots 1}
For $1 \leq k \leq n-1$, we have
\[
\zeta_i \cdot \alpha_k = \left\{
\begin{array}{r l}
\tilde{\alpha}_{k+i} & \quad \mbox{if} \ i \ \mbox{is even}, \\ 
\hat{\alpha}_{k+i} & \quad \mbox{if} \ i \ \mbox{is odd}, 
\end{array}
\right.
\] 
\enlem

\Proof
Using Lemma \ref{lem: aux root 1}, we observe that the set $\{ \tilde{\alpha}_k\}_{1 \leq k \leq n}$ is closed under the action of $s_0 [2,n{\shortminus}1][1,n{\shortminus}2]s_n$, and that this action sends $\tilde{\alpha}_k \mapsto \tilde{\alpha}_{k+1}$ (with indices modulo $n$). Repeating this action $\frac{i}{2}$-times, when $i$ is even, yields the first case of the lemma. 

Moreover, it follows that, if $i$ is odd, then $\big(s_0[2,n{\shortminus}1][1,n{\shortminus}2]s_n\big)^{\frac{i{\shortminus}1}{2}}$ acts by sending $\tilde{\alpha}_k \mapsto \tilde{\alpha}_{k+i-1}$. 
Using Lemma \ref{lem: aux root 1} again, we observe that $\pi_1[1,n{\shortminus}2]s_n$ defines a bijection 
\[
\{ \tilde{\alpha}_k\}_{1 \leq k \leq n} \to \{ \hat{\alpha}_k\}_{1 \leq k \leq n}, \qquad \tilde{\alpha}_k \mapsto \hat{\alpha}_{k+1}. \qedhere
\]
\enproof 

\subsubsection{Weights $\omega_{n-1}, \omega_n$}

Set 
\[
\zeta_i = \left\{
\begin{array}{r l}
\pi_{n{\shortminus}1} r_{n{\shortminus}2}r_{n{\shortminus}4} \cdots r_3 s_n 
& \quad \mbox{if }\  i = n{\shortminus}1 \ \mbox{is even}, \\ 
\pi_{n{\shortminus}1} s_{n{\shortminus}1}s_{n{\shortminus}2} r_{n{\shortminus}3} r_{n{\shortminus}5} \cdots r_3 s_n
& \quad \mbox{if }\  i = n{\shortminus}1 \ \mbox{is odd}, \\ 
\pi_{n} r_{n{\shortminus}2} r_{n{\shortminus}4} \cdots r_4 s_n [2,n{\shortminus}1]  
& \quad \mbox{if }\  i = n \ \mbox{is even}, \\ 
\pi_{n} s_{n{\shortminus}1}s_{n{\shortminus}2}r_{n{\shortminus}3} r_{n{\shortminus}5} \cdots r_4s_n [2,n{\shortminus}1]
& \quad \mbox{if }\  i = n \ \mbox{is odd}, \\ 
\end{array}
\right.
\]

\Lem \label{lem: action on sroots n,n-1}
We have 
\begin{align} 
\zeta_{n {\shortminus} 1} \cdot \alpha_k =& \ \left\{
\begin{array}{r l}
\alpha_{k{\shortminus}1} & \quad \mbox{if} \ 2 \leq k \leq n{\shortminus}1, \\ 
\pi_{n{\shortminus}1} [n {\shortminus} 1, 2] \cdot \alpha_{1} & \quad \mbox{if} \ k=1 \ \mbox{is odd}, 
\end{array}
\right. \\
\zeta_{n} \cdot \alpha_k =& \ \left\{
\begin{array}{r l}
\alpha_{k{\shortminus}1} & \quad \mbox{if} \ 2 \leq k \leq n{\shortminus}2, \\ 
\alpha_{n{\shortminus}2} & \quad \mbox{if} \ k=n, \\ 
\pi_{n} [n {\shortminus} 1, 2] \cdot \alpha_{1} & \quad \mbox{if} \ k=1 \ \mbox{is odd}. 
\end{array}
\right.
\end{align} 
\enlem

\Proof
The lemma follows by a direct calculation.  
\enproof

\subsection{Type $\mathsf{B}_n^{(1)}$} 

For $1 \leq i \leq n$, set 
\begin{align} \label{eq: zeta B}
\zeta_i = \begin{cases} 
(s_0[2,n][1,n])^{\frac{i}{2}}, & \text{if $i$ is even}, \\
\pi_1 [1,n](s_0[2,n][1,n])^{\frac{i-1}{2}} & \text{if $i$ is odd}. 
\end{cases}
\end{align}
Set 
\begin{align*}
\tilde{\alpha}_k = \begin{cases}
\alpha_k & 1 \leq k \leq n-1, \\
\alpha_0 + 2\alpha_n + \sum_{j=2}^{n-1} \alpha_j & k=0. 
\end{cases}
\end{align*}

\begin{lem} \label{lem: action on sroots 1 B}
If $1 \leq k \leq n-1$ then $\zeta_i.\alpha_k = \tilde{\alpha}_{i+k}$, where the indices are taken modulo~$n$.
\end{lem}

\Proof
The proof is similar to that of Lemma \ref{lem: action on sroots n,n-1}. 
\enproof

\subsection{Type $\mathsf{C}_n^{(1)}$} 

Let $\zeta_i = (s_0s_1\cdots s_n)^i$, for $1 \leq i \leq n-1$, and set  
\begin{align*}
\tilde{\alpha}_k = \begin{cases}
\alpha_k & 1 \leq k \leq n-1, \\
\sum_{j=0}^{n} \alpha_j & k=0. 
\end{cases}
\end{align*}

\begin{lem} \label{lem: action on sroots 1 C} 
For $1 \leq i,k \leq n-1$, we have $\zeta_i.\alpha_k = \tilde{\alpha}_{i+k}$, where the indices are taken modulo $n$.
\end{lem}

\Proof
The lemma follows by a direct calculation.  
\enproof


\nc{\te}{\widetilde{E}}
\nc{\qi}{q_i}
\nc{\tE}{\widetilde{E}}

\section{Good polynomials}
\label{sec: good polynomials}

The goal of this section is to formulate explicit criteria which will allow us to deduce the compatibility of the new Drinfeld presentations for $\Uu$ and $\Ui$. 
From now on assume that $\mathbf{s} = (0, \cdots, 0)$ and abbreviate $\widetilde{E}_i = - q_i^{-2} \KK_i E_i K_i\mi$. 

\subsection{Iterated $q$-brackets} 
\label{subsec: iterated w br}

As shown in \S \ref{sec-comb}, reduced expressions for the fundamental weights in all classical types contain subexpressions almost identical to the formulae encountered in type $\mathsf{A}$. Below we describe the braid group action corresponding to such subexpressions in terms of iterated $q$-brackets. 

Define non-commutative polynomials $P_k(y_1, \cdots, y_k)$ and $P'_k(y_1, \cdots, y_k)$ over $\C$ by induction in the following way: 
\begin{alignat}{3}
P_1(y_1) =& \ y_1, \qquad& P_{k+1}(y_1, \cdots, y_{k+1}) =& \ P_k(y_1, \cdots, y_{k-1}, [y_k,y_{k+1}]_q), \\
P'_1(y_1) =& \ y_1, \qquad& P'_{k+1}(y_1, \cdots, y_{k+1}) =& \ [P'_k(y_1, \cdots, y_k), y_{k+1}]_q.
\end{alignat}
Clearly, we have 
\begin{align}
P_k(y_1, \cdots, y_k) =& \ P_{l+1}(y_1, \cdots, y_l, P_{k-l}(y_{l+1}, \cdots, y_k)), \\
P'_k(y_1, \cdots, y_k) =& \ P'_{k-l+1}(P'_l(y_1, \cdots, y_l), y_{l+1}, \cdots, y_k). 
\end{align}
for any $1 \leq l \leq k-1$. 

We say that a tuple $(y_1, \cdots, y_k)$ is \emph{almost commuting} if $y_my_n = y_ny_m$, for $|m-n| > 1$.

\Lem \label{lem: PnP'}
The polynomials $P_k$ and $P'_k$ are equal if $(y_1, \cdots, y_k)$ is almost commuting. In particular, in that case, 
\[
P_2(y_1, P_2(y_2, y_3)) = P'_2(P'_2(y_1, y_2), y_3).
\]
\enlem

\Proof 
The first statement is proven by induction. By definition, $P_1 = P'_1$ and $P_2 = P'_2$. Moreover, 
\begin{align}
P_{k+1}(y_1, \cdots, y_{k+1}) =& \ P_2(y_1, P_k(y_2, \cdots,y_{k+1})) \\ 
=& \  P_2(y_1, P'_k(y_2, \cdots,y_{k+1})) \\ 
=& \ P'_k(P_2(y_1, y_2), y_3, \cdots, y_{k+1}) \\
=& \  P'_k(P'_2(y_1, y_2), y_3, \cdots, y_{k+1}) 
= P'_{k+1}(y_1, \cdots, y_{k+1}),
\end{align}
where in the second equality we used induction, and in the third equality the fact that the tuple is almost commuting. 
For the second part, observe that 
\[ 
P_2(y_1, P_2(y_2, y_3)) = P_3(y_1, y_2, y_3) 
= P'_3(y_1, y_2, y_3) = P'_2(P'_2(y_1,y_2), y_3). 
\qedhere
\]  
\enproof

\Lem \label{lem: exchange}
If $y_my_{m+1} = y_{m+1}y_m$ then
\[
P_k(\cdots, y_m, y_{m+1}, \cdots) = P_k(\cdots, y_{m+1}, y_{m}, \cdots). 
\]
\enlem 

\Proof
It suffices to consider the case $k=3$ with $[y_1, y_2] = 0$. Then
\begin{align}
P_3(y_1, y_2, y_3) =& \ y_1y_2y_3 - q(y_1y_3y_2 + y_2y_3y_1) + q^2 y_3 y_2 y_1 \\
=& \ y_2y_1y_3 - q(y_1y_3y_2 + y_2y_3y_1) + q^2 y_3 y_1 y_2 = P_3(y_2, y_1, y_3).  \qedhere
\end{align}
\enproof

We say that an ordered tuple $(i_1, \cdots, i_k)$ of mutually distinct elements of $\indx$ is a \emph{chain} if $a_{i_l,i_{l+1}} a_{i_{l+1},i_{l}}  =  1$ for each $1 \leq l \leq k{\shortminus}1$. 
We can identify the Dynkin subdiagram formed by $i_1, \cdots, i_k$ with a Dynkin diagram of type $\mathsf{A}_k$. To avoid double subscripts, we relabel $i_1 = 1$, etc. 
Let us abbreviate 
\eq \label{eq: tau defin}
\tau_{l} = 
{[{k{\shortminus}l+1}, k]} \cdots [1,l]s_l = {[{k{\shortminus}l+1}, k]} \cdots [2,l+1][1,{l{\shortminus}1}],
\eneq
for $1 \leq l \leq k$. 

\Lem \label{lem: type A chains}
Let $(1, \cdots, k)$ be a chain in $\indx$. 
\be
\item We have
\[
\Tbr_{[k,2]} (B_{1}) = P_k(B_{1}, \cdots, B_{k}). 
\]
Similarly, 
\[
T_{[k,2]} (F_{1}) = P_k(F_{1}, \cdots, F_{k}), \qquad 
T_{[k,2]} (\tE_{1}) = P_k(\tE_{1}, \cdots, \tE_{k}). 
\]
\item For $1 \leq l \leq k$, we have
\[
\Tbr_{\tau_{l}} (B_{l}) = P_l(B_{k}, \cdots, B_{k-l+2}, P_{k-l+1}(B_{1}, \cdots, B_{k-l+1})). 
\]
The same formula holds if $\Tbr_{\tau_{l}}$ is replaced by $T_{\tau_{l}}$, and every $B_i$ uniformly replaced by $F_i$ or $\tE_i$. 
\ee
\enlem 

\Proof
This is proven in \cite[Lemmas 9.3--9.4, Proposition 9.6]{Przez-23}. 
\enproof

Later, we will also need the polynomials
\eq \label{eq: bold P pols}
\mathbf{P}_i(a,b) = \sum_{r=0}^2 q^r b^{(2-r)} a b^{(r)} + \KK_i a, \qquad 
\widehat{\mathbf{P}}(a,b) = \sum_{r=0}^2 q^r b^{(2-r)} a b^{(r)}, 
\eneq
for $i \in \indx$. 

\subsection{Good polynomials} 
\label{subsec: good pols}

Let $P(B_0, \cdots, B_n)$ be any non-commutative polynomial in the variables $B_0, \cdots, B_n$ with coefficients in $\C(\KK_i)_{i \in \indx}$. Considering this polynomial as an element of $\Uu$, and substituting for an occurrence of a letter $B_j$ either $\tE_j$ or $F_j$, we can write it as a sum of polynomials which are homogeneous\footnote{
For example, suppose that $P(B_0, \cdots, B_n) = B_2^2B_1$. Then 
\[ 
B_2^2B_1 = \tE^2_2\tE_1 + \tE_2^2 F_1 + (\tE_2F_2\tE_1 + F_2 \tE_2 \tE_1) 
+( \tE_2F_2F_1 + F_2 \tE_2 F_1) + F_2^2\tE_1 + F_2^2F_1,
\] 
i.e., the monomial $B_2^2B_1$ becomes a sum of six homogeneous polynomials. 
} 
in each of the variables $\tE_j$ and $F_j$: 
\eq \label{eq: subterms}
\eta(P(B_0, \cdots, B_n)) = \sum_{\underline{d}} P_{\underline{d}}(\tE_0, F_0, \cdots, \tE_n, F_n), 
\eneq 
where $\underline{d} = (d_j^\pm)_{j \in \indx}$, with $d_j^+ = \deg_{\tE_j}P_{\underline{d}}(\cdot)$ and $d_j^- = \deg_{F_j} P_{\underline{d}}(\cdot)$.\footnote{E.g., for $P_{\underline{d}}(\cdot) = \tE_2^2F_1$, one has $d_2^+ = 2$, $d_1^- = 1$.}
We call the summands on the RHS \emph{subterms}, and the corresponding tuple $\underline{d}$ the type of a subterm. We say that that a subterm is \emph{mixed} if there exist $k,l \in \indx$ such that $d_k^+, d_l^- > 0$; otherwise a subterm is called \emph{pure}. 

We will now introduce some technical assumptions and definitions with view to proving Proposition \ref{pro: good -> in U+}. 
From this point on let us assume that each subterm satisfies $d_0^+ + d_0^- = 1$. 
Let $P_+$ be the sum of all subterms with $d_0^+ = 1$, and $P_-$ be the sum of all subterms with $d_0^- = 1$. Moreover, let $P_{--}$ be the sum of all subterms which have maximal total degree and $d_j^+ = 0$ for all $j \in \indx$.

\Defi \label{defi: i-good}
Let $i \in \indx_0$. A polynomial $P(B_0, \cdots, B_n)$ is called: 
\be 
\item 
\emph{good} if every mixed subterm of $P(B_0, \cdots, B_n)$ with $d_0^+ > 0$ vanishes; 
\item 
$i$-\emph{good} if it satisfies the following conditions: 
\be
\item every subterm $P_{\underline{d}}(\tE_0, F_0, \cdots, \tE_n, F_n)$ satisfies\footnote{The inequality below should be read as a family of inequalities of the coefficients along each simple root.}: 
\eq \label{i good defi properties}
\sum_{1 \leq j \leq n} (d_j^+ + d_j^-) \alpha_{j} \leq \theta - \alpha_i, \qquad  
d_0^+ + d_0^- = 1, 
\eneq
\item $P_+ = T_{\omega'_i}(\tE_i)$, \ $P_{--} = T_{\omega'_i}(F_i)$. 
\ee 
\ee
\enDefi 

The importance of $i$-good polynomials is conveyed by the following proposition. 

\Pro \label{pro: good -> in U+}
Suppose that $\Tbr_{\omega'_i}(B_i)$ can be expressed as an $i$-good polynomial. 
Then 
\eq \label{eq: general propo TT mod}
\eta(\mathbf{T}_{\omega_i'}(B_i)) \equiv T_{\omega_i'}(\eta(B_i)) \quad \mod \Uu_{d_i \geq 1,+}, 
\eneq
and the diagram 
\eq \label{eq: comm diagram compatible}
\begin{tikzcd}
\Ui(\widehat{\mathfrak{sl}}_2) \arrow[r, "\eta", hookrightarrow] \arrow[d, "\iota_i", "\wr"'] & \widetilde{U}_q(L\mathfrak{sl}_2) \arrow[d, "\iota_i", "\wr"']   \\
\Ui_{[i]} \arrow[d, hookrightarrow]  & \Uu_{[i]} \arrow[d, hookrightarrow]  \\
\Ui \arrow[r, "\eta", hookrightarrow] & \Uu
\end{tikzcd}
\eneq
commutes modulo $\Uu_{\neq i,+}$. 
\enpro

\Proof 
Suppose that $\Tbr_{\omega'_i}(B_i) = P(B_0, \cdots, B_n)$ is an $i$-good polynomial. Definition \ref{defi: i-good} implies that 
\eq \label{eq: two sums} 
\eta(P(B_0, \cdots, B_n)) = P_+ + P_-  = T_{\omega'_i}(\tE_i) + P_-. 
\eneq 
Let us analyse $P_-$ further. It can be written as the sum of (i) $P_{--} = T_{\omega'_i}(F_i)$, (ii) pure subterms with $\sum_{1 \leq j \leq n} d_j^- \alpha_{j} < \theta - \alpha_i$, and (iii) mixed subterms, with $d_0^- = 1$ in each case. It follows from Lemma \ref{lem: deg of hr} and \eqref{i good defi properties} that the degrees of these subterms (in the Drinfeld gradation) are given by
\[
(i) \ \alpha_i , \quad (ii) \ \theta - \sum_{1 \leq j \leq n} d_j^- \alpha_{j} > \alpha_i, \quad (iii) \ \theta + \sum_{1 \leq j \leq n} (d_j^+ - d_j^-) \alpha_{j} > \alpha_i,
\]
in the respective cases. In particular, the subterms of kinds (ii) and (iii) lie in $\Uu_{d_i \geq 1,+}$. This completes the proof of \eqref{eq: general propo TT mod}. 

Let $\mu$ denote the composition of the eastward and southward arrows, and $\nu$ the composition of the southward and eastward arrows in \eqref{eq: comm diagram compatible}. We have $\mu(B_0) = \nu(B_0)$ and, by the first part of the proposition, $\mu(B_{-1}) \equiv \nu(B_{-1})$ modulo $U_{d_i=1,+}$. Since $B_0$ and $B_{-1}$ generate $\Ui(\widehat{\mathfrak{sl}}_2) $, and $\mu(B_0)$ contains terms of degree $\pm \alpha_i$, it follows that the diagram commutes modulo $\Uu_{\neq i,+}$. 
\enproof 


\section{Auxiliary calculations - braid group actions} 
\label{sec: aux cal brgr}

In section \S \ref{sec: weak com} we will express $\Tbr_{\omega'_i}(B_i)$ in terms of certain polynomials. 
The following calculations will be used to prove that those polynomials are $i$-good. 

\Lem \label{lem: FE vanishing}
If $a_{ji} = -1$, then: 
\be
\item $[F_j,\tE_i]_{q_i} = 0$, 
\item $T_i(F_j) =  [F_j, F_i]_{q_i}$ and $T_i(\tE_j) =  [\tE_j, \tE_i]_{q_i}$. 
\ee 
\enlem

\Proof
The first part follows from the fact that, by \eqref{eq: Ke eK rel}, 
\eq 
F_jE_iK_i\mi - q_iE_iK_i\mi F_j = q_iE_iK_i\mi F_j - q_iE_iK_i\mi F_j = 0.
\eneq 
In the second part, $T_i(F_j) =  [F_j, F_i]_{q_i}$ follows immediately from the definition of the braid group action, and 
\begin{align}
T_i(\tE_j) =& \ - q_i^{-2} T_i (\KK_j E_j K_j\mi) \\
=& \ q_i^{-3} \KK_i\KK_j (E_jE_i - q_iE_iE_j) K_i\mi K_j\mi \\ 
=& \ q_i^{-4} \KK_i\KK_j (E_jK_j\mi E_iK_i\mi  - q_iE_iK_i\mi E_jK_j\mi) \\
=& \  [\tE_j, \tE_i]_{q_i}. \qedhere
\end{align}
\enproof 

\Lem \label{lem: TjTiBj formula}
If $a_{ij}a_{ji} = 2$ then 
\[
\Tbr_j \Tbr_i (B_j) = 
\left\{
\begin{array}{r l}
\KK_i B_j + \sum_{r=0}^2 (-q)^r B_i^{(2-r)} B_j B_i^{(r)} & \quad \mbox{if } \ a_{ji} = -2, \\ \\
{[B_i, B_j]}_{q^2} & \quad \mbox{if }\ a_{ji} = -1. 
\end{array}
\right. 
\]
\enlem 

\Proof
Let $a_{ji} = -2$. Then $q_i=q$, $q_j=q^2$ and
\begin{align*}
{\bf T}_i(B_j) = [2]_q^{-1} B_jB_i^2 - qB_iB_jB_i+ q^2 [2]_q^{-1} B_i^2B_j + \mathbb{K}_i B_j. 
\end{align*}
Moreover, $T_j( \mathbb{K}_j B_i ) = \mathbb{K}_j B_i$ and
\begin{align*}
{\bf T}_j( [2]_q^{-1} B_jB_i^2) = & \ [2]_q^{-1} \mathbb{K}_j^{-1} B_jB_iB_jB_iB_j- q^2 [2]_q^{-1} \mathbb{K}_j^{-1}B_jB_iB_j^2B_i \\
& - q^2 [2]_q^{-1} \mathbb{K}_j^{-1} B_j^2B_i^2B_j + q^4 [2]_q^{-1} \mathbb{K}_j^{-1} B_j^2B_iB_jB_i, \\
 {\bf T}_j( - qB_iB_jB_i ) = & \ - q \mathbb{K}_j^{-1} B_iB_j^2B_iB_j + q^3 \mathbb{K}_j^{-1} B_iB_j^3B_i  \\
& + q^3 \mathbb{K}_j^{-1} B_jB_iB_jB_iB_j- q^5 \mathbb{K}_j^{-1} B_jB_iB_j^2B_i, \\
{\bf T}_j( q^2 [2]_q^{-1} B_i^2B_j) =  & \ q^2  [2]_q^{-1} \mathbb{K}_j^{-1} B_iB_jB_iB_j^2-  q^4  [2]_q^{-1} \mathbb{K}_j^{-1} B_iB_j^2B_iB_j \\
& -  q^4  [2]_q^{-1}  \mathbb{K}_j^{-1} B_jB_i^2B_j^2 + q^6  [2]_q^{-1} \mathbb{K}_j^{-1} B_jB_iB_jB_iB_j.
\end{align*}
Therefore, 
\begin{align*}
{\bf T}_j{\bf T}_i(B_j) = & \ q \left(q^4 + 1\right) \mathbb{K}_j^{-1} B_jB_iB_jB_iB_j- q^2 \left(q^4 + q^2 + 1\right) [2]_q^{-1} \mathbb{K}_j^{-1} B_jB_iB_j^2B_i \\
& - q^2 [2]_q^{-1}  \mathbb{K}_j^{-1} B_j^2B_i^2B_j + q^4 [2]_q^{-1} \mathbb{K}_j^{-1} B_j^2B_iB_jB_i \\
&  - \left(q^4 + q^2 + 1\right) [2]_q^{-1} \mathbb{K}_j^{-1} B_iB_j^2B_iB_j + q^3 \mathbb{K}_j^{-1}B_iB_j^3B_i \\
& q^2 [2]_q^{-1} \mathbb{K}_j^{-1} B_iB_jB_iB_j^2 -  q^4 [2]_q^{-1} \mathbb{K}_j^{-1} B_jB_i^2B_j^2 + \mathbb{K}_i B_j.
\end{align*}
Repeatedly using the relation
\begin{align*}
B_j^2 B_i & = (q_j+q_j^{-1}) B_jB_iB_j - B_iB_j^2 - q_j^{-1} \mathbb{K}_j B_i \\
& = (q^2+q^{-2}) B_jB_iB_j - B_iB_j^2 - q^{-2} \mathbb{K}_j B_i, 
\end{align*}
we obtain 
\begin{align*}
& {\bf T}_j{\bf T}_i(B_j) = q^2 [2]_q^{-1} B_jB_i^2 + [2]_q^{-1} B_i^2B_j- qB_iB_jB_i + \mathbb{K}_iB_j.
\end{align*}

Next, let $a_{ji}=-1$. Then $q_i = q^2$ and $q_j = q$. We have ${\bf T}_i(B_j) = B_jB_i- q^2B_iB_j$ and
\begin{align*}
{\bf T}_j {\bf T}_i (B_j) = & \  [2]_q^{-1} \mathbb{K}_j^{-1} B_jB_iB_j^2 - q \mathbb{K}_j^{-1} B_j^2B_iB_j
+ 
  q^2 [2]_q^{-1} \mathbb{K}_j^{-1} B_j^3B_i+  B_jB_i \\
& - q^2 [2]_q^{-1} \mathbb{K}_j^{-1} B_i B_j^3 
 +q^3\mathbb{K}_j^{-1} B_jB_iB_j^2
- q^4 [2]_q^{-1} \mathbb{K}_j^{-1} B_j^2B_i  B_j- q^2 B_iB_j \\
 = & \ \left(q^4 + q^2 + 1\right) [2]_q^{-1} \mathbb{K}_j^{-1} B_jB_iB_j^2 - \left(q^4 + q^2 + 1\right) [2]_q^{-1} \mathbb{K}_j^{-1} B_j^2B_iB_j \\
& + q^2 [2]_q^{-1} \mathbb{K}_j^{-1} B_j^3B_i + B_jB_i- q^2  [2]_q^{-1} \mathbb{K}_j^{-1} B_iB_j^3- q^2B_iB_j.
\end{align*}
Using the Serre-type relation
\begin{align*}
B_j^3 B_i = & (q^4 + q^2 + 1) q^{-2} B_j^2B_iB_j-  \left(q^4 + q^2 + 1\right) q^{-2} B_jB_iB_j^2 + B_iB_j^3 \\
& - [2]_q^2 q^{-1} \mathbb{K}_jB_jB_i + [2]_q^2 q^{-1}  \mathbb{K}_jB_iB_j, 
\end{align*}
it follows that 
\begin{align*}
{\bf T}_j {\bf T}_i (B_j) = [B_i, B_j]_{q^2}, 
\end{align*}
completing the proof. 
\enproof 

\Lem \label{lem: TjTiBj appl} 
Let $a_{ij}a_{ji} = 2$. 
If $a_{ji} = -2$ then 
\begin{align}
\KK_i \tE_j + \sum_{r=0}^2 (-q)^r B_i^{(2-r)} \tE_j B_i^{(r)} = T_jT_i(\tE_j) =& \ \sum_{r=0}^2 (-q)^r \tE_i^{(2-r)}\tE_j \tE_i^{(r)}, \\
T_jT_i(F_j) =& \ \sum_{r=0}^2 (-q)^r F_i^{(2-r)} F_j F_i^{(r)}. 
\end{align}
If $a_{ji} = -1$ then 
\[
 [B_i, \tE_j]_{q^2} = T_jT_i(\tE_j) =  [\tE_i, \tE_j]_{q^2} , \qquad   T_jT_i(F_j) = [F_i, F_j]_{q^2}. 
\]
\enlem 

\Proof 
We prove the $a_{ji} = -2$ case, leaving the other easier case to the reader. 
We have $q_i = q$ and $q_j = q^2$.  
First, observe that 
\eq \label{eq: 1 phase BE calc}
\sum_{r=0}^2 (-q)^r F_i^{(2-r)} \tE_j F_i^{(r)} = ([2]\mi - q\mi +q^{-2}[2]\mi) F_i^2 \tE_j = 0. 
\eneq
Secondly, 
\begin{align}
[2]\mi(\tE_iF_i + F_i\tE_i)\tE_j - q(\tE_i \tE_j F_i + F_i\tE_j \tE_i) + q^2[2]\mi \tE_j(\tE_iF_i + F_i\tE_i) =& \\ 
[2]\mi(\tE_iF_i + F_i\tE_i)\tE_j - q\mi\tE_iF_i \tE_j - q^3\tE_j F_i \tE_i + q^2[2]\mi \tE_j(\tE_iF_i + F_i\tE_i) =& \\ 
-[2]\mi(q^{-2}\tE_i F_i - F_i\tE_i) \tE_j + [2]\mi \tE_j(q^2\tE_i F_i - q^4 F_i \tE_i) =& \\ 
-q^{-2} [2]\mi \KK_i\big( [E_i,F_i] K_i\mi\tE_j - q^4\tE_j [E_i,F_i] K_i\mi \big) =& \\ 
-q^{-2} (q^2 - q^{-2})\mi \KK_i \big( (1 - K_i^{-2})\tE_j - q^4\tE_j  (1 - K_i^{-2}) \big) =& \ -\KK_i \tE_j. \label{eq: 2 phase BE calc}
\end{align} 

Next, we calculate $T_jT_i(\tE_j)$. By definition, 
\begin{align*}
T_i(E_j) =  [2]\mi E_i^2E_j - q\mi E_iE_jE_i + q^{-2}[2]\mi E_jE_i^2,
\end{align*}
and 
\begin{align*}
T_jT_i(E_j) & = - [2]\mi E_jE_iE_jE_iF_jK_j + q^{-2}[2]\mi E_jE_i^2E_jF_jK_j \\
& + q^{-2}[2]\mi E_iE_j^2E_iF_jK_j  - q^{-4}[2]\mi E_iE_jE_iE_jF_jK_j \\ 
& + q\mi E_jE_iF_jK_jE_jE_i - q^{-3} E_jE_iF_jK_jE_iE_j - q^{-3}E_iE_jF_jK_jE_jE_i \\
& + q^{-5} E_iE_jF_jK_jE_iE_j - q^{-2}[2]\mi F_jK_jE_jE_iE_jE_i + q^{-4}[2]\mi F_jK_jE_jE_i^2E_j \\
& + q^{-4}[2]\mi F_jK_jE_iE_j^2E_i - q^{-6}[2]\mi F_jK_jE_iE_jE_iE_j. 
\end{align*}
Applying relations \eqref{eq: Ke eK rel}--\eqref{eq: Ei Ej rel}, the expression above simplifies to 
\begin{align*}
T_jT_i(E_j) = - q\mi E_iE_jE_i + [2]\mi E_jE_i^2 + q^{-2}[2]\mi E_i^2E_j. 
\end{align*}

Since $T_1T_2(\mathbb{K}_1) = \mathbb{K}_2^2 \mathbb{K}_1$ and $T_1T_2(K_1) =  K_2^2 K_1$, we obtain 
\begin{align} \label{eq: 3 phase BE calc}
\ \ \ T_jT_i(\tE_j) =& \  q^{-6} \sum_{r=0}^2 (-q)^r E_i^{(2-r)}E_j E_i^{(r)}K_i^{-2}K_j\mi \KK_i^2 \KK_j 
= \sum_{r=0}^2 (-q)^r \tE_i^{(2-r)}\tE_j \tE_i^{(r)}, 
\end{align}
where the second equality follows from the fact that $E_i^{(2-r)}E_j E_i^{(r)}K_i^{-2}K_j\mi = \linebreak   (E_iK_i\mi)^{(2-r)} E_jK_j\mi (E_iK_i\mi)^{(r)}$ for each $r$. 

Now, combining \eqref{eq: 1 phase BE calc}, \eqref{eq: 2 phase BE calc} and \eqref{eq: 3 phase BE calc}, we obtain the first equality in the statement of the lemma. The second equality is obtained via a similar calculation to that of $T_jT_i(E_j)$ above. 
\commentout{
In type $B_2^{(1)}$, let $i=2$ and $j=1$. Then $a_{ij}=-2$ and  $q_1=q^2$, $q_2=q$, and
\begin{align*}
T_2(F_1) =   \frac{q}{q^2 + 1}F_1F_2F_2- qF_2F_1F_2+\frac{q^3}{q^2 + 1}F_2F_2F_1,
\end{align*}
and
\begin{align*}
T_1T_2(F_1) & = - \frac{q}{\left(q^2 + 1\right)}\frac{1}{K_1}E_1F_2F_1F_2F_1+\frac{q^3}{q^2 + 1}\frac{1}{K_1}E_1F_2F_1F_1F_2+\frac{q^3}{q^2 + 1}\frac{1}{K_1}E_1F_1F_2F_2F_1 \\
& - \frac{q^5}{\left(q^2 + 1\right)}\frac{1}{K_1}E_1F_1F_2F_1F_2+qF_2F_1\frac{1}{K_1}E_1F_2F_1- q^3F_2F_1\frac{1}{K_1}E_1F_1F_2 \\
& - q^3F_1F_2\frac{1}{K_1}E_1F_2F_1+q^5F_1F_2\frac{1}{K_1}E_1F_1F_2 - \frac{q^3}{\left(q^2 + 1\right)}F_2F_1F_2F_1\frac{1}{K_1}E_1 \\
& +\frac{q^5}{q^2 + 1}F_2F_1F_1F_2\frac{1}{K_1}E_1 +\frac{q^5}{q^2 + 1}F_1F_2F_2F_1\frac{1}{K_1}E_1- \frac{q^7}{\left(q^2 + 1\right)}F_1F_2F_1F_2\frac{1}{K_1}E_1. 
\end{align*}
Applying relations \refeq{eq: Ke eK rel}--{eq: Ei Ej rel}, the expression above simplifies to 
\begin{align*}
K_i e_j^{\pm} = q_i^{\pm a_{ij}} e_j^{\pm} K_i, \quad [e_i^+, e_j^-] = \delta_{ij} \frac{K_i-K_i^{-1}}{q_i-q_i^{-1}},
\end{align*} 
we have that
\begin{align*}
T_1T_2(F_1) =  - qF_2F_1F_2+\frac{q}{q^2 + 1}F_2F_2F_1+\frac{q^3}{q^2 + 1}F_1F_2F_2. 
\end{align*}
}
\enproof

\section{Weak compatibility}
\label{sec: weak com}

The goal of this section is to prove the following theorem. 

\Thm \label{thm: weak comp}
In types $\mathsf{A}_n^{(1)}, \mathsf{B}_n^{(1)}, \mathsf{C}_n^{(1)}, \mathsf{D}_n^{(1)}$, for each $1 \leq i \leq n$, $\Tbr_{\omega'_i}(B_i)$ can be expressed as an $i$-good polynomial. Hence
\[
\eta(\mathbf{T}_{\omega_i'}(B_i)) \equiv T_{\omega_i'}(\eta(B_i)) \mod \Uu_{d_i \geq 1,+}. 
\]
\enthm

The second statement follows immediately from the first by Proposition \ref{pro: good -> in U+}. 
To prove the first statement, we proceed by a direct, case-by-case method. We first explicitly describe the polynomials expressing $\Tbr_{\omega'_i}(B_i)$, using the calculations from \S \ref{sec: root combin},
and then prove that these polynomials are $i$-good, using the calculations from \S \ref{sec: aux cal brgr}. 
Type $\mathsf{A}_n^{(1)}$ has already been handled in \cite[Proposition 9.6]{Przez-23}. 

\subsection{Type $\mathsf{D}_n^{(1)}$} 
For convenience, we will sometimes write $B'_l = B_l$ for $2 \leq l \leq n{\shortminus}2$ and $B'_{n{\shortminus}1} = B_n$. We use similar notation with $B$ replaced by~$F$. 

\subsubsection{Weights $\omega_1, \cdots, \omega_{n{\shortminus}2}$}
\label{para: D ord wts brgr}

\Lem \label{lem: T B pol D}
Let $1 \leq i \leq n{\shortminus}2$. 
If $i$ is even, then 
\eq \label{eq: D even pol}
\Tbr_{\omega'_i}(B_i) = 
P_i(B_{i{\shortminus}1}, \cdots, B_{1}, P_{n{\shortminus}i}(B_{i+1}, \cdots, B_{n{\shortminus}1}, P_{n{\shortminus}1}(B_{n},B_{n{\shortminus}2}, \cdots, B_2, B_0))). 
\eneq
If $i$ is odd, then 
\eq \label{eq: D odd pol}
\Tbr_{\omega'_i}(B_i) = 
P_i(B_{i{\shortminus}1}, \cdots, B_{1}, P_{n{\shortminus}i}(B_{i+1}, \cdots, B_{n{\shortminus}2}, B_{n}, P_{n{\shortminus}1}(B_{n{\shortminus}1}, \cdots, B_2, B_0))). 
\eneq
\enlem

\Proof 
Recall the definitions of $\zeta_i$ from \S \ref{sec: root combin} and $\tau_i$ from \eqref{eq: tau defin}. 
Since $\omega'_i = \zeta_i \tau_i$ and $\ell(\omega'_i) = \ell(\zeta_i) + \ell(\tau_i)$, 
we have $\Tbr_{\omega'_i}(B_i) = \Tbr_{\zeta_i} \Tbr_{\tau_{i}} (B_{i})$. 
According to Lemma \ref{lem: type A chains} (with $k=n{\shortminus}1$ and $l=i$), 
\eq \label{eq: T tau i}
\Tbr_{\tau_{i}} (B_{i}) = P_i(B_{n{\shortminus}1}, \cdots, B_{n{\shortminus}i+1}, P_{n{\shortminus}i}(B_{1}, \cdots, B_{n{\shortminus}i})). 
\eneq
Using Lemmas \ref{lem: br act on wt} and \ref{lem: action on sroots 1}, we observe that applying $\Tbr_{\zeta_i}$ to the RHS of \eqref{eq: T tau i} merely changes the indices of the variables except the last variable $B_{n{\shortminus}i}$. The precise change of the indices is given by Lemma \ref{lem: action on sroots 1}, and 
\[
\Tbr_{\zeta_i}(B_{n{\shortminus}i}) = \left\{ 
\begin{array}{r l} 
P_{n{\shortminus}1}(B_{n},B_{n{\shortminus}2}, \cdots, B_2, B_0) & \quad \mbox{if} \ i \ \mbox{is even}, \\ 
P_{n{\shortminus}1}(B_{n{\shortminus}1},B_{n{\shortminus}2}, \cdots, B_2, B_0) & \quad \mbox{if} \ i \ \mbox{is odd}, \\
\end{array}
\right.
\]
follows directly from the definition of the braid group action. 
\enproof

\Pro \label{pro: D i-good 1}
The polynomials \eqref{eq: D even pol}--\eqref{eq: D odd pol} are $i$-good. 
\enpro

\Proof
We do the even case, the odd one being analogous. 
First consider the innermost nested polynomial $P_{n{\shortminus}1}(B_{n},B_{n{\shortminus}2}, \cdots, B_2, B_0)$. Define 
\[
R_k = \left\{
\begin{array}{r l}
B_0 & \quad \mbox{if} \ k=1, \\
P_k(B_k, \cdots, B_2, B_0) & \quad \mbox{if} \ 2 \leq k \leq n{\shortminus}2, \\
P_{n{\shortminus}1}(B_n, B_{n{\shortminus}2}, \cdots, B_2, B_0) & \quad \mbox{if} \ k = n{\shortminus}1.
\end{array}
\right.
\]
Let $S_k$ be defined analogously, with each $B_l$ replaced by $\tE_l$. 
We will prove by induction that the polynomials $R_k$ are good. Observe that $R_1$ is good trivially, and $R_2$ is good by part (1) of Lemma \ref{lem: FE vanishing}. For $k \geq 2$, 
we have 
\[
R_{k+1} = P_{3}(B'_{k+1}, B_k, R_{k{\shortminus}1}). 
\]  
According to the inductive hypothesis, $R_{k{\shortminus}1}$ and $R_k$ are good. 
Therefore, it suffices to show that $P_2(F'_{k+1}, S_k) = 0$. 
Since $F'_{k+1}$ commutes with $S_{k{\shortminus}1}$, Lemma \ref{lem: PnP'} implies that 
\eq \label{eq: FS vanish}
P_2(F'_{k+1}, S_k) = P'_2([F'_{k+1}, \tE_k]_q, S_{k{\shortminus}1}), 
\eneq
which vanishes since $[F'_{k+1}, \tE_k]_q = 0$ by part (1) of Lemma \ref{lem: FE vanishing}. 
We conclude that $R_{k+1}$ is good and, hence, by induction, $R_{n{\shortminus}1}$ is good.

Next, we show that the middle nested polynomial $P_{n{\shortminus}i}(B_{i+1}, \cdots, B_{n{\shortminus}1}, R_{n{\shortminus}1})$ is good. 
The proof is similar to the argument above. 
Set
\[
R_k^{\mathsf{mid}} = \left\{
\begin{array}{r l}
R_{n{\shortminus}1} & \quad \mbox{if} \ k=0, \\
P_{k+1}(B_{n{\shortminus}k}, \cdots, B_{n{\shortminus}1}, R_{n{\shortminus}1}) & \quad \mbox{if} \ 1 \leq k \leq n{\shortminus}i{\shortminus}1,
\end{array}
\right.
\]
and define $S_k^{\mathsf{mid}}$ analogously, with each $B_l$ replaced by $\tE_l$, and $R_{n{\shortminus}1}$ by $S_{n{\shortminus}1}$. 
We already know that $R_0^{\mathsf{mid}}$ is good. The polynomial $R_1^{\mathsf{mid}}$ is also good since, by Lemma \ref{lem: exchange}, 
\[
P_2(F_{n{\shortminus}1}, S_{n{\shortminus}1}) = P_2(\tE_n, P_2(F_{n{\shortminus}1}, S_{n{\shortminus}2})), 
\]
and $P_2(F_{n{\shortminus}1}, S_{n{\shortminus}2}) = 0$ by the same argument as in \eqref{eq: FS vanish}. 

Next, we claim that $[F_{n{\shortminus}k}, S_{k{\shortminus}2}^{\mathsf{mid}}] = 0$ for each $2 \leq k \leq n{\shortminus}i{\shortminus}1$. Since $S_{k{\shortminus}2}^{\mathsf{mid}} = P_{k{\shortminus}1}(\tE_{n{\shortminus}k+2}, \cdots, \tE_{n{\shortminus}1}, S_{n{\shortminus}1})$, it suffices to show that $[F_{n{\shortminus}k}, S_{n{\shortminus}1}] = 0$. This is indeed true because 
\[
[F_{n{\shortminus}k}, S_{n{\shortminus}1}] = T_{\zeta_i} ( [F_{n{\shortminus}k{\shortminus}i}, \tE_{n{\shortminus}i}] ) = 0. 
\] 
Therefore, 
\[
P_2(F_{n{\shortminus}k}, S_{k{\shortminus}1}^{\mathsf{mid}}) = P'_2([F_{n{\shortminus}k}, \tE_{n{\shortminus}k+1}]_q, S_{k{\shortminus}2}^{\mathsf{mid}}) = 0, 
\]
which implies that $R_{k}^{\mathsf{mid}}$ is good. By induction, we conclude that $R_{n{\shortminus}i{\shortminus}1}^{\mathsf{mid}}$ is good.

Finally, we show that the outer polynomial $P_i(B_{i{\shortminus}1}, \cdots, B_{1}, R_{n{\shortminus}i{\shortminus}1}^{\mathsf{mid}})$ is good. 
Set
\[
R_k^{\mathsf{out}} = \left\{
\begin{array}{r l}
R_{n{\shortminus}i{\shortminus}1}^{\mathsf{mid}} & \quad \mbox{if} \ k=0, \\
P_{k+1}(B_{k}, \cdots, B_{1}, R_{n{\shortminus}i{\shortminus}1}^{\mathsf{mid}}) & \quad \mbox{if} \ 1 \leq k \leq i{\shortminus}1,
\end{array}
\right.
\]
and define $S_k^{\mathsf{out}}$ analogously, with each $B_l$ replaced by $\tE_l$, and $R_{n{\shortminus}i{\shortminus}1}^{\mathsf{mid}}$ by $S_{n{\shortminus}i{\shortminus}1}^{\mathsf{mid}}$. 
We already know that $R_0^{\mathsf{out}}$ is good. The polynomial $R_1^{\mathsf{out}}$ is also good since
\begin{align}
P_2(F_{1}, S_{n{\shortminus}i{\shortminus}1}^{\mathsf{mid}}) =& \ 
T_{\zeta_i}(P_2(F_{n{\shortminus}i+1}, P_{n{\shortminus}i}(\tE_1, \cdots, \tE_{n{\shortminus}i} ) )) \\
=& \ T_{\zeta_i}(P_{n{\shortminus}i}(\tE_1, \cdots, \tE_{n{\shortminus}i{\shortminus}1}, P_2(F_{n{\shortminus}i+1}, \tE_{n{\shortminus}i}))), 
\end{align} 
and $P_2(F_{n{\shortminus}i+1}, \tE_{n{\shortminus}i}) = 0$ by part (1) of Lemma \ref{lem: FE vanishing}. 
We also claim that $[F_{k}, S_{k{\shortminus}2}^{\mathsf{out}}] = 0$ for each $2 \leq k \leq i{\shortminus}1$. Since $S_{k{\shortminus}2}^{\mathsf{out}} = P_{k{\shortminus}1}(\tE_{k{\shortminus}2}, \cdots, \tE_{1}, S_{n{\shortminus}i{\shortminus}1}^{\mathsf{mid}})$, it suffices to show that $[F_{k}, S_{n{\shortminus}i{\shortminus}1}^{\mathsf{mid}}] = 0$. This is indeed true because 
\[
[F_{k}, S_{n{\shortminus}i{\shortminus}1}^{\mathsf{mid}}] = T_{\zeta_i} ( [F_{n+k{\shortminus}i}, P_{n{\shortminus}i}(\tE_1, \cdots, \tE_{n{\shortminus}i})] ) = 0. 
\]
It follows that  
\[
P_2(F_{k}, S_{k{\shortminus}1}^{\mathsf{out}}) = P'_2([F_{k}, \tE_{k{\shortminus}1}]_q, S_{k{\shortminus}2}^{\mathsf{out}}) = 0, 
\]
which implies that $R_{k}^{\mathsf{out}}$ is good. By induction, we conclude that $R_{i{\shortminus}1}^{\mathsf{out}}$ is good. 
This proves that \eqref{eq: D even pol} is a good polynomial.

Next, Lemma \ref{lem: br act on wt}, part (2) of Lemma \ref{lem: type A chains}, and part (2) of Lemma \ref{lem: FE vanishing} imply that, by uniformly substituting $F_j$ or $\tE_j$ for each $B_j$ in the proof of Lemma \ref{lem: T B pol D}, we obtain a proof of the following identities:
\begin{align}
T_{\omega'_i}(F_i) =& \
P_i(F_{i{\shortminus}1}, \cdots, F_{1}, P_{n{\shortminus}i}(F_{i+1}, \cdots, F_{n{\shortminus}1}, P_{n{\shortminus}1}(F_{n},F_{n{\shortminus}2}, \cdots, F_2, F_0))), \\
T_{\omega'_i}(\tE_i) =& \
P_i(\tE_{i{\shortminus}1}, \cdots, \tE_{1}, P_{n{\shortminus}i}(\tE_{i+1}, \cdots, \tE_{n{\shortminus}1}, P_{n{\shortminus}1}(\tE_{n},\tE_{n{\shortminus}2}, \cdots, \tE_2, \tE_0))). 
\end{align}

The fact that \eqref{eq: D even pol} is a good polynomial, together with the two identities above, yields part (2.b) of the 
Definition \ref{defi: i-good} of $i$-goodness. Part (2.a) clearly holds as well, which completes the proof. 
\enproof

\subsubsection{Weights $\omega_{n-1}, \omega_n$} 
\label{para: D sp wts brgr}

We apply similar analysis to weights $\omega_{n-1}$ and $\omega_n$. 

\Pro
We have 
\begin{align} 
\label{eq: Dn-1 pol}
\Tbr_{\omega'_{n \smin 1}}(B_{n \smin 1}) =& \
P_{n \smin 1}(B_{n \smin 2}, \cdots, B_1, P_{n \smin 1}(B_n, B_{n \smin 2}, \cdots, B_2, B_0)), \\
\label{eq: Dn pol}
\Tbr_{\omega'_n}(B_n) =& \
P_{n \smin 1}(B_{n \smin 2}, \cdots, B_1, P_{n \smin 1}(B_{n \smin 1}, B_{n \smin 2}, \cdots, B_2, B_0)). 
\end{align} 
The polynomials \eqref{eq: Dn-1 pol}--\eqref{eq: Dn pol} are $i$-good. 
\enpro

\Proof
Let $i \in \{n \smin 1, n\}$.  
Since $\omega'_i = \zeta_i [1,n\smin 2]$ and $\ell(\omega'_i) = \ell(\zeta_i) + \ell([1,n\smin 2])$, 
we have $\Tbr_{\omega'_i}(B_i) = \Tbr_{\zeta_i} \Tbr_{[1,n\smin 2]} (B_{i})$. 
According to Lemma \ref{lem: type A chains} (with $k=n{\shortminus}1$ and $l=i$), 
\eq \label{eq: T tau n,n-1}
\Tbr_{[1,n\smin 2]} (B_{i}) = P_{n\smin1}(B'_{n{\shortminus}1}, B_{n\smin2}, \cdots, B_1) 
\eneq 
Using Lemma \ref{lem: action on sroots n,n-1}, we observe that applying $\Tbr_{\zeta_i}$ to the RHS of \eqref{eq: T tau n,n-1} merely shifts the indices of the variables by $-1$, except the last variable $B_{1}$, and 
\[
\Tbr_{\zeta_i}(B_1) = \left\{ 
\begin{array}{r l} 
P_{n{\shortminus}1}(B_{n},B_{n{\shortminus}2}, \cdots, B_2, B_0) & \quad \mbox{if} \ i=n \smin 1, \\ 
P_{n{\shortminus}1}(B_{n{\shortminus}1},B_{n{\shortminus}2}, \cdots, B_2, B_0) & \quad \mbox{if} \ i=n. \\
\end{array}
\right.
\]
The proof of $i$-goodness follows the same method as the proof of Proposition \ref{pro: D i-good 1}, so we omit it. 
\enproof

\subsection{Type $\mathsf{B}_n^{(1)}$}

In this subsection only we consider polynomials $P_k(y_1, \cdots, y_k)$ to be defined in the same way as in \S \ref{subsec: iterated w br} but with $q$ replaced by $q^2$. We will also use the polynomials $\mathbf{P}_i(a,b)$ from \eqref{eq: bold P pols}.

\Lem \label{lem: T B pol B}
Let $1 \leq i \leq n{\shortminus}1$. Then 
\eq \label{eq: B good pol}
\Tbr_{\omega'_i}(B_i) = 
P_i(B_{i{\shortminus}1}, \cdots, B_{1}, P_{n{\shortminus}i}(B_{i+1}, \cdots, B_{n{\shortminus}1}, \mathbf{P}_n( P_{n{\shortminus}1}(B_{n{\shortminus}1}, \cdots, B_2, B_0), B_n))).
\eneq 
If $i=n$ then
\eq \label{eq: B good pol 2}
\Tbr_{\omega'_n}(B_n) = P_n(B_{n-1}, \cdots, B_1, P_n(B_n, \cdots, B_2, B_0)).  
\eneq
\enlem 

\Proof 
Let $1 \leq i \leq n{\shortminus}1$. 
Arguing in the same way as in the proof of Lemma \ref{lem: T B pol D}, we find that 
\eq \label{eq: T tau i B}
\Tbr_{\omega'_i}(B_i) = \Tbr_{\zeta_i}(P_i(B_{n{\shortminus}1}, \cdots, B_{n{\shortminus}i+1}, P_{n{\shortminus}i}(B_{1}, \cdots, B_{n{\shortminus}i}))), 
\eneq
with $\zeta_i$ given by \eqref{eq: zeta B}. Lemma \ref{lem: action on sroots 1 B} implies that $\Tbr_{\zeta_i}$ shifts the indices of all the variables on the RHS of \eqref{eq: T tau i B} by $i$ modulo $n$, except the last variable $B_{n{\shortminus}i}$. 
Moreover, 
\[
\Tbr_{\zeta_i}(B_{n{\shortminus}i}) = 
\left\{ 
\begin{array}{r l} 
\Tbr_{s_0[2,n]}(B_{n{\shortminus}1}) & \quad \mbox{if} \ i \ \mbox{is even}, \\ 
\Tbr_{\pi_1[1,n]}(B_{n{\shortminus}1}) & \quad \mbox{if} \ i \ \mbox{is odd}. 
\end{array}
\right. 
\] 
By Lemma \ref{lem: TjTiBj formula}, 
\begin{align}
\Tbr_{s_0[2,n]}(B_{n{\shortminus}1}) =& \ \Tbr_{s_0[2,n-2]}\mathbf{P}_n(B_{n{\shortminus}1}, B_n) \\ 
=& \ \mathbf{P}_n(P_{n{\shortminus}1}(B_{n{\shortminus}1}, \cdots, B_2, B_0), B_n)  
\end{align}
if $i$ is even, with an analogous calculation in the odd case. The proof for $i=n$ is similar to that in type $\mathsf{C}_n^{(1)}$, which is done in Lemma \ref{lem: T B pol C}, so we leave it to the reader. 
\enproof

\Pro
The polynomials \eqref{eq: B good pol}--\eqref{eq: B good pol 2} are $i$-good. 
\enpro 

\Proof 
The proof is similar to the proof of Proposition \ref{pro: D i-good 1}. Let us explain the differences in the case $1 \leq i \leq n - 1$. 
Arguing as in the first paragraph of the proof of Proposition \ref{pro: D i-good 1}, we deduce that the polynomials $P_{k}(B_{k}, \cdots, B_2, B_0)$, for $2 \leq k \leq n-1$, are good. 
Hence $S_k := P_{k}(B_{k}, \cdots, B_2, \tE_0) = P_{k}(\tE_{k}, \cdots, \tE_2, \tE_0)$. Then 
\begin{align}
\mathbf{P}_n(S_{n-1}, B_n) = \mathbf{P}_n(P_2(\tE_{n-1}, S_{n-2}), B_n) = P_2(\mathbf{P}_n(\tE_{n-1}, B_n), S_{n-2}), 
\end{align}
and $\mathbf{P}(\tE_{n-1}, B_n) = \widehat{\mathbf{P}}(\tE_{n-1},\tE_n)$ by Lemma \ref{lem: TjTiBj appl}. 
Therefore, 
\[
\mathbf{P}_n(S_{n-1}, B_n) = P_2(\widehat{\mathbf{P}}(\tE_{n-1},\tE_n), S_{n-2}) = \widehat{\mathbf{P}}(S_{n-1}, \tE_n). 
\]
The rest of the argument follows the same pattern as in the proof of \ref{pro: D i-good 1}. In particular, we obtain 
\begin{align}
 P_i(B_{i{\shortminus}1}, \cdots, B_{1}, P_{n{\shortminus}i}(B_{i+1}, \cdots, B_{n{\shortminus}1}, \mathbf{P}_n( P_{n{\shortminus}1}(B_{n{\shortminus}1}, \cdots, B_2, \tE_0), B_n))) =& \ \\
 P_i(\tE_{i{\shortminus}1}, \cdots, \tE_{1}, P_{n{\shortminus}i}(\tE_{i+1}, \cdots, \tE_{n{\shortminus}1}, \widehat{\mathbf{P}}( P_{n{\shortminus}1}(\tE_{n{\shortminus}1}, \cdots, \tE_2, \tE_0), \tE_n))) =& \ T_{\omega'_i}(\tE_i), \\ \label{eq: F T(F) type B proof} 
\ P_i(F_{i{\shortminus}1}, \cdots, F_{1}, P_{n{\shortminus}i}(F_{i+1}, \cdots, F_{n{\shortminus}1}, \widehat{\mathbf{P}}( P_{n{\shortminus}1}(F_{n{\shortminus}1}, \cdots, F_2, F_0), F_n))) =& \ T_{\omega'_i}(F_i). 
\end{align} 
Note that \eqref{eq: F T(F) type B proof} is indeed a subterm of \eqref{eq: B good pol} of maximal total degree, and so $P_{--} = T_{\omega'_i}(F_i)$, 
proving that \eqref{eq: B good pol} is $i$-good. 
\enproof

\subsection{Type $\mathsf{C}_n^{(1)}$} 

Finally, let us consider type $\mathsf{C}_n^{(1)}$. 

\Lem \label{lem: T B pol C}
Let $1 \leq i \leq n - 1$. Then 
\eq \label{eq: C good pol}
\Tbr_{\omega'_i}(B_i) = 
P_i(B_{i{\shortminus}1}, \cdots, B_{1}, P_{n{\shortminus}i}(B_{i+1}, \cdots, B_{n{\shortminus}1}, \mathbb{P}(B_n, P_{n}(B_{n{\shortminus}1}, \cdots, B_2, B_1, B_0)))), 
\eneq
where $\mathbb{P}(a,b) = [a,b]_{q^2}$. 
If $i=n$ then
\eq \label{eq: C good pol 2}
\Tbr_{\omega'_n}(B_n) =  \mathbf{P}_{n-1}( \cdots \mathbf{P}_2(\mathbf{P}_1(B_0, B_1), B_2), \cdots, B_{n-1}). 
\eneq
\enlem

\Proof
In the first case, the proof is the same as the proof of Lemma \ref{lem: T B pol B}, with the exception that Lemma \ref{lem: action on sroots 1 C} is used instead of Lemma \ref{lem: action on sroots 1 B}, and the case of $a_{ji} = -1$ (rather than $a_{ji} = -2$)  from Lemma \ref{lem: TjTiBj formula} is used. 

If $i = n$, then one proceeds as follows. One can rewrite \eqref{eq: om n C} as
\[
\omega_n = \pi_n [n,1] \cdots [n,n-1] s_n. 
\]
By Lemma \ref{lem: TjTiBj formula}, $\mathbf{T}_{[n,n-1]}(B_n) = \mathbf{P}_{n-1}(B_n,B_{n-1})$. Proceeding by induction, we may assume that 
\[
\mathbf{T}_{[n,k] \cdots [n,n-1]}(B_n) = \mathbf{P}_{k}( \cdots \mathbf{P}_{n-2}(\mathbf{P}_{n-1}(B_n, B_{n-1}), B_{n-2}), \cdots, B_{k}). 
\]
By Lemma \ref{lem: br act on wt}, applying $\mathbf{T}_{[n,k-1]}$ shifts the indices of all the $B_{l}$ and $\mathbf{P}_l$ by $-1$, except for $B_n$, which is transformed into $\mathbf{P}_{n-1}(B_n,B_{n-1})$. Finally, $\mathbf{T}_{\pi_n}$ reverses the indices. 
\enproof

\Pro
The polynomials \eqref{eq: C good pol}--\eqref{eq: C good pol 2} are $i$-good. 
\enpro 

\Proof
The proof in the case $1 \leq i \leq n - 1$ is similar to those in types $\mathsf{B}_n^{(1)}$ and $\mathsf{D}_n^{(1)}$, so let us consider the case $i = n$. First, observe that an appropriate modification of the proof of Lemma \ref{lem: T B pol C} yields 
\begin{align} 
{T}_{[n,k] \cdots [n,n-1]}(\tE_n) =& \ \widehat{\mathbf{P}}( \cdots \widehat{\mathbf{P}}(\widehat{\mathbf{P}}(\tE_n, \tE_{n-1}), \tE_{n-2}), \cdots, \tE_{k}) =: S_k, \\ \label{eq: C omega n P hat formula E}
T_{\omega'_n}(\tE_n) =& \ \widehat{\mathbf{P}}( \cdots \widehat{\mathbf{P}}(\widehat{\mathbf{P}}(\tE_0, \tE_1), \tE_2), \cdots, \tE_{n-1}), \\ \label{eq: C omega n P hat formula F}
T_{\omega'_n}(F_n) =& \ \widehat{\mathbf{P}}( \cdots \widehat{\mathbf{P}}(\widehat{\mathbf{P}}(F_0, F_1), F_2), \cdots, F_{n-1}), 
\end{align} 
for $1 \leq k \leq n-1$. 

We prove by (descending) induction that 
\eq \label{eq: bold P ind}
R_k := \mathbf{P}_{k}( \cdots \mathbf{P}_{n-2}(\mathbf{P}_{n-1}(\tE_n, B_{n-1}), B_{n-2}), \cdots, B_{k}) = S_k. 
\eneq
The base case of $k = n-1$ reduces to Lemma \ref{lem: TjTiBj appl}. Next, suppose that \eqref{eq: bold P ind} holds for~$k$. Then 
\begin{align}
R_{k-1} = \mathbf{P}_{k-1}(R_k, B_{k-1}) =& \ \mathbf{P}_{k-1}(S_k, B_{k-1}) \\ 
=& \ \KK_{k-1} S_k + \sum_{r=0}^2 (-q)^r B_{k-1}^{(2-r)} S_k B_{k-1}^{(r)}. 
\end{align}

Let us compute the sum above.  
Observe that $F_{k-1} S_k = q^2 S_k F_{k-1}$ (since $F_{k-1}$ commutes with all the letters in $S_k$ except for $K_k\mi$, which appears twice). Hence 
\begin{align}
\sum_{r=0}^2 (-q)^r F_{k-1}^{(2-r)} S_k F_{k-1}^{(r)} =& \ [2]\mi (q^4 - (q+q\mi)q^3 +q^2)S_k F_{k-1}^2 = 0. 
\end{align} 
Moreover, an easy calculation shows that 
\begin{align}
(F_{k-1}\tE_{k-1} + \tE_{k-1} F_{k-1}) S_k 
\ -& \ q(F_{k-1} S_k \tE_{k-1} + \tE_{k-1} S_k F_{k-1})  \\ 
+& \ q^2 S_k (F_{k-1} \tE_{k-1} + \tE_{k-1} F_{k-1}) 
= - \KK_{k-1} S_k. 
\end{align} 
It follows that $\mathbf{P}_{k-1}(S_k, B_{k-1}) = \widehat{\mathbf{P}}(S_k, \tE_{k-1}) = S_{k-1}$, concluding the inductive step. 

The induction above yields $R_1 = S_1$. Applying $\mathbf{T}_{\pi_n}$, together with \eqref{eq: C omega n P hat formula E}, yields $P_+ = T_{\omega'_i}(\tE_i)$, while \eqref{eq: C omega n P hat formula F} yields $P_{--} = T_{\omega'_i}(F_i)$. It follows that \eqref{eq: C good pol 2} is $i$-good. 
\enproof


\section{Strong compatibility}
\label{sec: strong compat}

The goal of this section is to prove a generalization of Theorem \ref{thm: rank 1 factorization} to the classical types $\mathsf{B}_n^{(1)}, \mathsf{C}_n^{(1)}, \mathsf{D}_n^{(1)}$. 

\Thm \label{thm: ultimate factorization theorem} 
Let $\widehat{\g}$ be of type $\mathsf{A}_n^{(1)}, \mathsf{B}_n^{(1)}, \mathsf{C}_n^{(1)}$ or~$\mathsf{D}_n^{(1)}$. Then 
\begin{align}
\eta_{\mathbf{s}}(\Thgsr_i(z)) \equiv& \ \xi_{\mathbf{s}}(\Thgsr_i(z)) \pmb{\phi}_i^-(z\mi)\pmb{\phi}_i^+(\CCC z) \fext{\mod \Uu_{+}}{z}, \label{eq: coproduct formula strong thm 1} \\
\Delta_{\mathbf{s}}(\Thgsr_i(z)) \equiv& \ \eta_{\mathbf{s}}(\Thgsr_i(z)) \otimes \eta(\Thgsr_i(z)) \quad  \mod \fext{\Uu  \otimes \Uu_{+}}{z}.  \label{eq: coproduct formula strong thm 2}
\end{align} 
\enthm

Type $\mathsf{A}_n^{(1)}$ was already handled in \cite[Corollary 9.15]{Przez-23}. Here we will follow the same methodology as in \emph{op.\ cit.} We begin in \S \ref{subsec: gen fact} below by recalling results which generalize to our situation without any extra work.

\subsection{Generalized factorization}
\label{subsec: gen fact}

In this subsection let $\widehat{\g}$ be of type $\mathsf{A}_n^{(1)}, \mathsf{B}_n^{(1)}, \mathsf{C}_n^{(1)}$ or~$\mathsf{D}_n^{(1)}$. 
Set $\CCC_i = \CCC\mi \mathbb{K}_i$. 

\Lem \label{lem: str comp basic}
For each $i \in \indx_0$: 
\[
\eta(A_{i,-1}) = x^-_{i,1} - \CCC_ix^+_{i,-1}K_i + Q_i, \qquad Q_i \in \Uu_{d_i \geq 1,+}. 
\]
Hence 
\[
\eta(A_{i,-1}) \equiv \iota_i \eta(A_{-1}) \quad \mod  \Uu_{d_i \geq 1,+}. 
\]
Moreover, 
\[
\eta(H_{i,1}) \equiv \iota_i\eta(H_1) + q_i^2 \CCC_i\mi [Q_i, F_i]_{q_i^{-2}} \quad \mod \Uu_{d_i \geq 2,+}. 
\]
\enlem 

\Proof
See \cite[Lemmas 9.8--9.9]{Przez-23}. 
\enproof 

The statement of Lemma \ref{lem: str comp basic} should also be seen as the definition of the term $Q_i$. 
As in \S \ref{subsec: good pols}, let us write  $\Tbr_{\omega'_i}(B_i)  = P_+ + P_-$. 
It then follows from the proof of \cite[Lemma 9.8]{Przez-23} that, explicitly, 
\eq \label{eq: Qi expl}
Q_i = \CCC_i (P_{-} - P_{--}). 
\eneq 

\Pro
Let $r\in \Z$. If 
\eq \label{eq: condition main}
[[Q_i, F_i]_{q_i^{-2}}, x_{i,-r}^-] \in \Uu_{d_i \geq 1,+},  
\eneq 
then 
\[
\eta(A_{i,r}) \equiv \iota_i \eta(A_{r}) \quad \mod  \Uu_{d_i \geq 1,+}. 
\]
\enpro 

\Proof 
See \cite[Proposition 9.10]{Przez-23}. 
\enproof

\Lem \label{cor: ultimate factorization theorem} 
If \eqref{eq: condition main} is true, then Theorem \ref{thm: ultimate factorization theorem} holds. 
\enlem

\Proof 
The proof is the same as in \cite[Corollary 9.15]{Przez-23}. 
\enproof

Therefore, it suffices to prove that condition \eqref{eq: condition main} holds in each type. 

\subsection{Explicit computation} 

Below we calculate explicitly that \eqref{eq: condition main} holds in type $\mathsf{D}_n^{(1)}$, for weights $\omega_{1}, \cdots, \omega_{n-2}$. The other types and weights can be handled using the same methods. 

Given \eqref{eq: Qi expl}, an explicit expression for $Q_i$ can be obtained from \eqref{eq: D even pol}--\eqref{eq: D odd pol}\footnote{For simplicity, we use formula \eqref{eq: D even pol}. In the odd case, indices $n$ and $n-1$ should be switched.}. Substituting either $F_j$ or $\tE_j$ for each occurrence of a variable $B_j$, we can write $Q_i$ as a sum of polynomials in the variables $F_j$ and $\tE_j$. In the series of lemmas below, we check that condition \eqref{eq: condition main} holds for each of these polynomials. Depending on the type of polynomial $P$, we prove one of the following: (i) $P$ vanishes, (ii) $[P, F_i]_{q^{-2}} \in \Uu_{d_i \geq 2,+}$, (iii) $P$ can, modulo $\Uu_{d_i \geq 2,+}$, be expressed as a polynomials in $x_{j,r}^+$ $(j \neq i)$. 

Let $\tilde{e}_i^+ = \tE_i$ and $\tilde{e}_i^- = F_i$. 
Let us also abbreviate 
\[
R = P_{n{\shortminus}1}(B_{n},B_{n{\shortminus}2}, \cdots, B_2, F_0). 
\]

\Lem
If $2 \leq i \leq n-2$, then   
\[ 
[P_i(\tE_{i{\shortminus}1}, B_{i{\shortminus}2}, \cdots, B_{1}, P_{n{\shortminus}i}(\tE_{i+1}, B_{i+2} \cdots, B_{n{\shortminus}1}, R)), F_i]_{q^{-2}}  \in  \Uu_{d_i \geq 2,+}. 
\]
\enlem

\Proof 
We calculate  
\begin{align} 
[P_i(\tE_{i{\shortminus}1}, B_{i{\shortminus}2}, \cdots, B_{1}, P_{n{\shortminus}i}(\tE_{i+1}, B_{i+2} \cdots, B_{n{\shortminus}1}, R)), F_i]_{q^{-2}} =& \\
[\tE_{i-1}, P_{i-1}(B_{i-2}, \cdots, B_1, [\tE_{i+1}, P_{n-i-1}(B_{i+1}, \cdots, B_{n-1}, [F_i,R])] ].& 
\end{align}  
For degree reasons, we only need to check that the commutator of $F_i$ with 
\[
R' = \sum_{\substack{\varepsilon_k \le \varepsilon_{k+1} \\ \in\{\pm\}}} P_{n{\shortminus}1}(\tilde{e}_{n}^{\varepsilon_n}, \tilde{e}_{n-2}^{\varepsilon_{n-2}}, \cdots, \tilde{e}_{i+1}^{\varepsilon_{i+1}},F_i, \cdots, F_2, F_0) 
\]
vanishes. This reduces to an easily verifiable calculation in the finite quantum group of type $\mathsf{A}_3$. 
\enproof

\Lem
If $1 \leq i \leq n-2$, then 
\begin{align} \label{eq: hard1}
P_i(B_{i{\shortminus}1}, \cdots, B_{1}, P_{n{\shortminus}i}(F_{i+1}, B_{i+2}, \cdots, B_{n{\shortminus}1}, P_{n{\shortminus}1}(B_{n},B_{n{\shortminus}2}, \cdots, B_2, F_0)))&
= \\
\sum_{\varepsilon_k \le \varepsilon_{k+1}} P_i(\tilde{e}_{i-1}^{\varepsilon_{i-1}}, \cdots, \tilde{e}_{1}^{\varepsilon_1}, P_{n{\shortminus}i}(F_{i+1}, \cdots, F_{n{\shortminus}1}, P_{n{\shortminus}1}(F_{n},F_{n{\shortminus}2}, \cdots, F_2, F_0)))& \\   \label{eq: hard2}
+ \delta_{i,n-2} \sum_{\varepsilon_k \le \varepsilon_{k+1}} P_i(\tilde{e}_{n-3}^{\varepsilon_{i-3}}, \cdots, \tilde{e}_{1}^{\varepsilon_1}, P_{n{\shortminus}i}(\tE_{n}, P_{n{\shortminus}1}(F_{n-1},F_{n{\shortminus}2}, \cdots, F_2, F_0)))&. 
\end{align}
\enlem 

\Proof 
By \cite[Lemmas 9.3--9.4]{Przez-23}, the LHS of \eqref{eq: hard1} reduces to 
\[
\sum_{\substack{\varepsilon_k \le \varepsilon_{k+1}, \\ \varepsilon'_k \le \varepsilon'_{k+1}}} P_i(\tilde{e}_{i-1}^{\varepsilon'_{i-1}}, \cdots, \tilde{e}_{1}^{\varepsilon'_1}, P_{n{\shortminus}i}(F_{i+1}, F_{i+2}, \cdots, F_{n{\shortminus}1}, P_{n{\shortminus}1}(\tilde{e}_{n}^{\varepsilon_n}, \tilde{e}_{n-2}^{\varepsilon_{n-2}}, \cdots, \tilde{e}_{2}^{\varepsilon_2}, F_0))).
\]
The lemma now follows from the following two equalities 
\begin{align} \label{eq: hard3}
[F_{n-1}, P_{n{\shortminus}1}(\tE_{n}, \tE_{n{\shortminus}2}, B_{n{\shortminus}3}, \cdots, B_2, F_0) ]_q =& \ 0, \\
\label{eq: hard4}
[F_{n-2}, [F_{n-1}, P_{n{\shortminus}1}(\tE_{n}, F_{n{\shortminus}2}, B_{n{\shortminus}3}, \cdots, B_2, F_0) ]_q]_q =& \ 0, 
\end{align} 
which can be checked directly\footnote{Such equalities can easily be verified using, e.g., GAP's QuaGroup Package \cite{GAP}. For an explicit code, see \href{https://github.com/lijr07/identities\_in\_quantum\_groups}{\sf https://github.com/lijr07/identities\_in\_quantum\_groups}. } (they reduce to calculations in finite quantum groups of type $\mathsf{A}_3$ and $\mathsf{D}_4$, respectively). 
\enproof 

\Lem
Let $2 \leq i \leq n-2$ and 
\begin{align}
Z =& \  \sum_{\varepsilon_k \le \varepsilon_{k+1}} P_i(\tE_{i-1}, \tilde{e}_{i-2}^{\varepsilon_{i-2}}, \cdots, \tilde{e}_{1}^{\varepsilon_1}, P_{n{\shortminus}i}(F_{i+1}, \cdots, F_{n{\shortminus}1}, P_{n{\shortminus}1}(F_{n},F_{n{\shortminus}2}, \cdots, F_2, F_0))). 
\end{align}
The element $[Z, F_i]_{q^{-2}}$ can be expressed as a polynomial in $x_{j,r}^+$ $(j \neq i)$, with coefficients in $\C$. Hence $[[Z, F_i]_{q^{-2}}, x_{i,s}^-] =0$,  for all $s \in \Z$. 
\enlem 

\Proof
The proof is analogous to \cite[Lemmas 9.13]{Przez-23}. 
\enproof

Let 
\[
Z' = \ \sum_{\substack{\varepsilon_k \le \varepsilon_{k+1}, \\ \varepsilon'_k \ge \varepsilon'_{k+1}}} P_i(F_{i-1}, \cdots, F_{1}, P_{n{\shortminus}i}(\tE_{i+1}, \tilde{e}_{i+2}^{\varepsilon'_{i+2}}, \cdots, \tilde{e}_{n-1}^{\varepsilon'_{n-1}}, P_{n{\shortminus}1}(\tilde{e}_{n}^{\varepsilon_n}, \tilde{e}_{n-2}^{\varepsilon_{n-2}}, \cdots, \tilde{e}_{2}^{\varepsilon_2}, F_0))).
\]

\Lem
Let $1 \leq i \leq n-2$. Then 
\begin{align}
Z' =& \ \sum_{\substack{ \varepsilon'_k \ge \varepsilon'_{k+1}}} P_i(F_{i-1}, \cdots, F_{1}, P_{n{\shortminus}i}(\tE_{i+1}, \tilde{e}_{i+2}^{\varepsilon'_{i+2}}, \cdots, \tilde{e}_{n-1}^{\varepsilon'_{n-1}}, P_{n{\shortminus}1}(F_{n}, F_{n-2}, \cdots, F_0))) \\ 
+& \ \sum_{\substack{\varepsilon_k \le \varepsilon_{k+1}}} P_i(F_{i-1}, \cdots, F_{1}, P_{n{\shortminus}i}(\tE_{i+1},\cdots, \tE_{n-1}, P_{n{\shortminus}1}(\tE_{n}, \tE_{n-2}, \tilde{e}_{n-3}^{\varepsilon_{n-3}}, \cdots, \tilde{e}_{2}^{\varepsilon_2}, F_0))) \\
+& \ P_i(F_{i-1}, \cdots, F_{1}, P_{n{\shortminus}i}(\tE_{i+1},\cdots, \tE_{n-2}, \tE_{n-1}, P_{n{\shortminus}1}(\tE_{n}, F_{n-2}, \cdots, F_2, F_0))) \\
+& \ \delta_{i \neq n-2} P_i(F_{i-1}, \cdots, F_{1}, P_{n{\shortminus}i}(\tE_{i+1},\cdots, \tE_{n-2}, \tE_{n}, P_{n{\shortminus}1}(F_{n-1}, F_{n-2}, \cdots, F_2, F_0))). 
\end{align}
\enlem

\Proof
This follows from \eqref{eq: hard3}--\eqref{eq: hard4}. 
\enproof

\Lem
The element $[Z', F_i]_{q^{-2}}$ can be expressed as a polynomials in $x_{j,r}^+$ $(j \neq i)$, with coefficients in $\C$, modulo $ \Uu_{d_i \geq 2,+}$. Hence $[[Z', F_i]_{q^{-2}}, x_{i,s}^-] = 0$ modulo $\Uu_{d_i \geq 1,+}$,  for all $s \in \Z$. The same holds if $Z'$ is replaced by \eqref{eq: hard2}. 
\enlem

\Proof
The proof is analogous to \cite[Lemmas 9.13]{Przez-23}. 
\enproof


\section{Application to $q$-characters}

We propose a generalization of the notion of $q$-characters to affine quantum symmetric pairs, based on the Lu--Wang Drinfeld-type presentation. We apply Theorem \ref{thm: ultimate factorization theorem} to show that our construction is compatible with the usual $q$-character map for quantum affine algebras. 

\subsection{$q$-characters of quantum affine algebras}

Fix $\mathbf{c} = (c_0, \cdots, c_n) \in (\C^\times)^{n+1}$ and let $C$ be the image of $\mathfrak{C}$ in $\Uic$. 
Let $\Rep \Uq$ be the monoidal category of finite dimensional representations of $\Uq$. It acts on $\Rep \Uic$, the category of finite dimensional representations of $\Uic$, via the coproduct $\Delta_{\mathbf{c},\mathbf{s}}$ \eqref{eq: coproduct explicit}. It follows from \cite[Sec.\ 3.5]{dobson-kolb-19} (see also \cite[Lemmas 2.4]{Przez-23}) that any such monoidal action is, up to twisting by a character, equal to the ``standard'' monoidal action with $\mathbf{s} = (0, \cdots, 0)$. Therefore, henceforth we shall assume that we are working with the standard monoidal action. 

We will use the same notations for the Grothendieck groups of the categories above. From this point of view, $\Rep \Uq$ is a ring acting on $\Rep \Uic$. 
Let $\Uh$ be the subalgebra of $\Uq$ generated by $h_{i,k}$ ($i \in \indx_0$, $k < 0$). Following \cite{FrenRes, FrenMuk-comb}, set 
\begin{align} \label{eq: defi of Yia}
Y_{i,a} =& \  K_{\omega_i}^{-1} \exp \left( -(q-q\mi) \sum_{k > 0} \tilde{h}_{i,-k} a^k z^k  \right) \in \Uhz \qquad (a \in \C^\times),
\end{align}
where 
\eq
\tilde{h}_{i,-k} = \sum_{j \in \indx_0} \widetilde{C}_{ji}(q^k) \hg_{j,-k},
\eneq 
and $\widetilde{C}(q)$ is the inverse of the $q$-Cartan matrix. 
Let $\Yring = \Z[Y_{i,a}^{\pm 1}]_{i \in \indx_0, a \in \C^{\times}}$. 
By \cite[Theorem 3]{FrenRes}, there exists an injective ring homomorphism $\chmap \colon \Rep \Uq \to \Yring \subset \Uhz$, called the \emph{$q$-character map}, given by 
\[
[V] \mapsto \Tr_V \left[ \exp \left( -(q-q\mi) \sum_{i \in \indx_0} \sum_{k >0} \frac{k}{[k]_{q_i}} \pi_V(h_{i,k}) \otimes \tilde{h}_{i,-k} z^k \right) \cdot (\pi_V \otimes 1)(T) \right], 
\]
where $T$ is as in \cite[(3.8)]{FrenRes}, and $\pi_V \colon \Uq \to \End(V)$ is the representation. 

The expression above derives from the Khoroshkin-Tolstoy-Levendorsky-Soibelman-Stukopkin-Damiani (KTLSSD) factorization of the universal $R$-matrix \cite{khor-tol, lev-sob-st, damiani-r}. However, by \cite[Proposition 2.4]{FrenMuk-comb}, the $q$-character map can equivalently be defined more explicitly in terms of the joint spectrum of the Drinfeld-Cartan operators (i.e., the coefficients of the series $\pmb{\phi}^{\pm}_i(z)$). More precisely, there is a one-to-one correspondence between the monomials occurring in $\chmap(V)$ and the common eigenvalues of $\pmb{\phi}^{\pm}_i(z)$ on $V$. 
Let 
$V = \bigoplus_{\gamma} V_{\gamma}$, with $\gamma = (\gamma^\pm_{i,\pm m})_{i \in \indx_0, m \in \Z_{\ge 0}}$, where 
\begin{align*}
V_{\gamma} = \{ v \in V \mid \exists p \ \forall i \in \indx_0 \ \forall m \in \Z_{\ge 0}:  (\psi^\pm_{i,\pm m} - \gamma^\pm_{i,\pm m})^p \cdot v = 0 \}, 
\end{align*}
be the generalized eigenspace decomposition of $V$. 
Collect the eigenvalues into generating series $\gamma^\pm_i(z) = \sum_{m \ge 0} \gamma^\pm_{i,\pm m} z^{\pm m}$. 
By \cite[Proposition 2.4]{FrenMuk-comb}, the series $\gamma^\pm_i(z)$ are expansions (at $0$ and $\infty$, respectively) of  the same rational function of the form
\[
q_i^{\deg Q_i - \deg R_i} \frac{Q_i(q_i^{-1}z)R_i(q_iz)}{Q_i(q_iz)R_i(q_i^{-1}z)}, 
\]
for some polynomials $Q_i(z), R_i(z)$ with constant term $1$. Writing  
\eq \label{eq: QR notation}
Q_i(z) = \prod_{r=1}^{k_i} (1 - z a_{i,r}), \qquad R_i(z) = \prod_{s=1}^{l_i} (1 - z b_{i,s}), 
\eneq
the $q$-character of $V$ can now be expressed as 
\[
\chi_q(V) = \sum_{\gamma} \dim(V_\gamma) M_\gamma, \qquad 
M_\gamma = \prod_{i \in \indx_0} \prod_{r=1}^{k_i} Y_{i,a_{i,r}} \prod_{s=1}^{l_i} Y_{i,b_{i,s}}\mi. 
\]

\subsection{Boundary $q$-characters} 

The notion of a quasi-$K$-matrix appeared first in the work of Bao and Wang \cite{BW18a, BW18b} on canonical bases for quantum symmetric pairs, as an intertwiner between the bar involutions on the quantum group and the coideal subalgebra.  This construction of the quasi-$K$-matrix was later extended to arbitrary Kac--Moody type by Balagovi\'c
and Kolb in \cite{BK19} and, in finite type, led to the realization of the universal $K$-matrix as a coideal
intertwiner \cite{BW18b, BK19}. 
Later, Appel and Vlaar \cite{AppelVlaar, AppelVlaar3} generalized these methods to obtain a universal $K$-matrix in the affine case. However, a 
KTLSSD-type factorization is not yet available.\footnote{In the finite case, a factorization of quasi $K$-matrices was established in \cite{dobson-kolb-19, WangZhang}.} Therefore, we propose to initiate the study of $q$-characters for affine quantum symmetric pairs based on the Lu--Wang (Drinfeld-type) presentation.

\Defi
Let  
\[
\mathcal{K}^0 = \frac{1-q_i^{-2}\CCC z^2}{1-\CCC z^2}  \exp \left( -(q-q\mi) \sum_{i \in \indx_0} \sum_{k >0} \frac{k}{[k]_{q_i}} H_{i,k} \otimes \tilde{h}_{i,-k} z^k \right) \in \fext{\Uic \otimes \Uh}{z}. 
\]
We define the \emph{boundary $q$-character map} to be 
\[
\ichmap \colon \Rep \Uic \ \to \ \Uhz, \qquad [V] \mapsto \Tr_V(\mathcal{K}^0 \circ (\pi_V \otimes 1)). 
\]
\enDefi

Consider $\Uhz$ as a $\Yring$-module via the ring homomorphism
\eq \label{eq: Yring module}
\Yring \to \Yring \hookrightarrow \Uhz, \qquad Y_{i,a} \mapsto Y_{i,Ca}Y_{i,a\mi}\mi. 
\eneq
Below we apply Theorem \ref{thm: ultimate factorization theorem} to show that our boundary $q$-character map is a module homomorphism. First, we need some notation. Given a polynomial $P(z) \in \C[z]$ with constant term $1$, let $P^\dag(z)$ be the polynomial with constant term $1$ whose roots are obtained from those of $P(z)$ via the transformation $a \mapsto \CC\mi a\mi$; and let $P^*(z)$ be the polynomial with constant term $1$ whose roots are the inverses of the roots of $P(z)$.

\Pro 
\label{cor: FR thm Oq} 
Let $W \in \Rep \Uq$. 
Then the generalized eigenvalues of $\pmb{\grave{\Thg}}_i(z)$ on the restricted representation $\eta_{\mathbf{c}}^*(W)$   are of the form: 
\eq \label{eq: FR formula Oq}
\gamma_i^\iota(z) = \frac{\ourQ_i(q_i^{-1}z)}{\ourQ_i(q_i z)} \frac{\ourQ_i^\dag(q_i z)}{\ourQ_i^\dag(q_i^{-1}z)}, 
\eneq
where $\ourQ_i(z)$ is a polynomial with constant term $1$. Explicitly, 
\[
\ourQ_i(z) = {Q_i(\CC z)R_i^*(z)}, \quad  \ourQ_i^\dag(z) = {R_i(\CC z)Q_i^*(z)}. 
\]
\enpro

\Proof 
Theorem \ref{thm: ultimate factorization theorem} implies that the action of $\pmb{\grave{\Thg}}_i(z)$ on $\eta_{\mathbf{c}}^*(W)$ is the sum of the action of $\pmb{\phi}_i^-(z\mi)\pmb{\phi}_i^+(C z)$ and a nilpotent operator. Hence the eigenvalues of $\pmb{\grave{\Thg}}_i(z)$ coincide with those of $\pmb{\phi}_i^-(z\mi)\pmb{\phi}_i^+(C z)$. The result now follows from \cite[Proposition 1]{FrenRes} by the same argument as in the proof of \cite[Corollary 5.1]{Przez-23}. 
\enproof

\Cor \label{cor: comm diagram qchar actions}
Let $(\Uq, \Uic)$ be a split quantum symmetric pair of type $\mathsf{A}_n^{(1)}, \mathsf{B}_n^{(1)}, \mathsf{C}_n^{(1)}$ or~$\mathsf{D}_n^{(1)}$. Then the following diagram commutes: 
\[
\begin{tikzcd}[ row sep = 0.2cm]
\Rep \Uq  \arrow[r, "\chmap"] & \Yring  \\
 \curvearrowright  & \curvearrowright  \\
\Rep \Uic \arrow[r, "\ichmap"] & \Uhz. 
\end{tikzcd}
\]
\encor

\Proof 
The proof is a modification of the proof of \cite[Proposition 2.4]{FrenMuk-comb}. 
Let $W \in \Rep \Uq$ and $V \in \Rep \Uic$. We need to calculate $\ichmap(V \otimes W)$. 
Formulae \eqref{eq: coproduct formula strong thm 1}--\eqref{eq: coproduct formula strong thm 2} in Theorem \ref{thm: ultimate factorization theorem}, together with Proposition \ref{cor: FR thm Oq}, imply that the eigenvalues of $\pmb{\grave{\Thg}}_i(z)$ on $V \otimes W$ are a product of the eigenvalues of $\pmb{\grave{\Thg}}_i(z)$ on $V$, and $\gamma_i^\iota(z)$. 
Hence $\Tr_{V\otimes W}(\mathcal{K}^0 \circ (\pi_{V \otimes W} \otimes 1)) = \Tr_V(\mathcal{K}^0 \circ (\pi_V \otimes 1)) \cdot \Tr_W(\mathcal{K}^0 \circ (\pi_W \otimes 1))$, and the eigenvalues of $H_{i,m}$ on $W$ are (up to the overall normalization factor $\frac{1-q_i^{-2}\CC z^2}{1-\CC z^2}$) of the form 
\eq \label{eq: eig with twist}
\frac{q_i^m - q_i^{-m}}{n(q-q\mi)}\left( \sum_{r=1}^{k_i} ((\CC a_{i,r})^m - a_{i,r}^{-m}) -
\sum_{s=1}^{l_i} ((\CC b_{i,s})^m - b_{i,s}^{-m}) \right),
\eneq
keeping the notation from \eqref{eq: QR notation}. Plugging \eqref{eq: eig with twist} into the definition of $\ichmap(W)$ and comparing with \eqref{eq: defi of Yia}, we conclude that $\ichmap(W)$ is a linear combination of monomials of the form 
\[
\prod_{i \in \indx_0} \prod_{r=1}^{k_i} Y_{i,C a_{i,r}} Y_{i,a_{i,r}}\mi \prod_{s=1}^{l_i} Y_{i, b_{i,s}} Y_{i,C b_{i,s}}\mi. 
\] 
Therefore, $\ichmap(W)$ is equal to the image of $\chmap(W)$ under \eqref{eq: Yring module}, completing the proof. 
\enproof

\Rem
In the case of the ordinary $q$-character map $\chi_q$, the property of being a ring homomorphism follows almost directly from its definition via the universal $R$-matrix (and the associated properties of transfer matrices). In contrast, in our case, the fact that $\ichmap$ intertwines the actions of $\Rep \Uq$ and $\Yring$ is non-trivial, since we are not using a universal $K$-matrix, but working with the Drinfeld-type presentation instead. 
\enrem

Finally, let us consider the rank one example of the $q$-Onsager algebra. 

\begin{exam} 
We calculate the boundary $q$-characters of restricted evaluation representations, for $(s_0,s_1) = (0,0)$.  
Let $W_n(a)$ be the irreducible $n+1$-dimensional representation of $U_q(\mathfrak{sl}_2)$ evaluated at $qa$. Then, by \cite[{(4.3)}]{FrenRes}, its $q$-character is 
\[
\chi_q(W_n(a)) = \sum_{i=0}^n M_i, \qquad M_i = \prod_{k=i+1}^n Y_{aq^{n-2k+1}} \prod_{k=1}^i Y_{aq^{n-2k+3}}\mi. 
\]
By Corollary \ref{cor: comm diagram qchar actions}, 
\[
\ichmap(W_n(a)) = \sum_{i=0}^n \mathbf{M}_i, \quad \mathbf{M}_i = \prod_{k=i+1}^n Y_{Caq^{n-2k+1}}Y_{a\mi q^{-n+2k-1}}\mi \prod_{k=1}^i Y_{a\mi q^{-n+2k-3}} Y_{Caq^{n-2k+3}}\mi. 
\]
Writing $\chi_q(W_n(q^{-2}C\mi a\mi)) = \sum_{i=0}^n M'_i$ and $\ichmap(W_n(q^{-2}C\mi a\mi)) = \sum_{i=0}^n \mathbf{M}'_i$, it is easy to see that $\mathbf{M}_i = \mathbf{M}'_{n-i}$, which implies that 
\[
\ichmap(W_n(a)) = \ichmap(W_n(q^{-2}C\mi a\mi)). 
\] 
This is not a coincidence since, by \cite[Theorem 1.17]{ito-ter-10}, $W_n(a)$ and $W_n(q^{-2}C\mi a\mi)$ are indeed isomorphic as representations of the $q$-Onsager algebra. 
\end{exam}

\addtocontents{toc}{\SkipTocEntry}
\section*{Conflict of interest statement} 

The authors have no competing interests to declare that are relevant to the content of this
article. 

\addtocontents{toc}{\SkipTocEntry}
\section*{Data availability statement} 

The authors declare that all the data supporting the findings of this article are available within the paper. 

\providecommand{\bysame}{\leavevmode\hbox to3em{\hrulefill}\thinspace}
\providecommand{\MR}{\relax\ifhmode\unskip\space\fi MR }
\providecommand{\MRhref}[2]{%
  \href{http://www.ams.org/mathscinet-getitem?mr=#1}{#2}
}
\providecommand{\href}[2]{#2}

\end{document}